\documentclass[11pt,leqno]{article}

\usepackage{graphicx}
\usepackage{latexsym} 
\usepackage{amsmath}
\usepackage{amssymb}
\usepackage{amsfonts}
\usepackage{enumerate}
\usepackage{theorem}

\theoremstyle{plain}
\newtheorem{teo}{Theorem}[section]
\newtheorem{lem}[teo]{Lemma}
\newtheorem{cor}[teo]{Corollary}
\newtheorem{prop}[teo]{Proposition}

{\theorembodyfont{\rmfamily}

\newtheorem{oss}[teo]{Remark}

\newtheorem{claim}{Claim}}

\renewcommand{\eqref}[1]{\textnormal{(\ref{#1})}}

\numberwithin{equation}{section}

\newcommand{\cvd}{\hfill$\square$}
\newcommand{\proof}[1]{\noindent\textsc{Proof#1}}
\newcommand{\rmi}{\mathrm{i}}
\newcommand{\rme}{\mathrm{e}}
\newcommand{\rmd}{\mathrm{d}}

\title{Analysis of an Inverse Problem Arising in Photolithography}

\author{Luca Rondi\thanks{Dipartimento di Matematica e Informatica,
Universit\`a degli Studi di Trieste, via Valerio, 12/1, 34127
Trieste, Italy. \texttt{rondi@units.it}} \and\
Fadil Santosa\thanks{School of Mathematics, University of Minnesota, 
Minneapolis, MN 55455, USA. \texttt{santosa@math.umn.edu} } }

\date{}

\begin{document}

\maketitle

\setcounter{section}{0}
\setcounter{secnumdepth}{2}

\begin{abstract}
We consider the inverse problem of determining an optical mask that
produces a desired circuit pattern in photolithography. We
set the problem as a shape design problem in which the unknown
is a two-dimensional domain.  The relationship between the target
shape and the unknown is modeled through diffractive optics.  We
develop a variational formulation that is well-posed and
propose an approximation that can be shown to have convergence
properties.  The approximate problem can serve as a
foundation to numerical methods.

\medskip


\medskip

\end{abstract}

\section{Introduction}\label{intro}

Photolithography is a key process in the production of integrated
circuits. It is the process by which circuit patterns are transferred
onto silicon wafers. A review of this manufacturing technology is given in
\cite{Schell03}.  The main step in photolithography is the
creation of a circuit image on the photoresist coating which sits on
the silicon layer that is to be patterned.  The image is formed using
ultra-violet (UV) light which is diffracted by a mask, and refracted by a
system of lenses.  The mask simply consists of cut-outs, and lets
light through the holes.  The parts of the photoresist that are exposed to
the UV light can be removed, leaving openings to the layer to be
patterned.  The next stage is etching, which removes material in the
layer that is unprotected by the photoresist.  Once etching is done,
the photoresist can be removed, and the etched away ``channels'' may
be filled.  The entire process is illustrated schematically in Figure 1.

\begin{figure}[t]
\hspace{2in}{\includegraphics[width=3.5in]{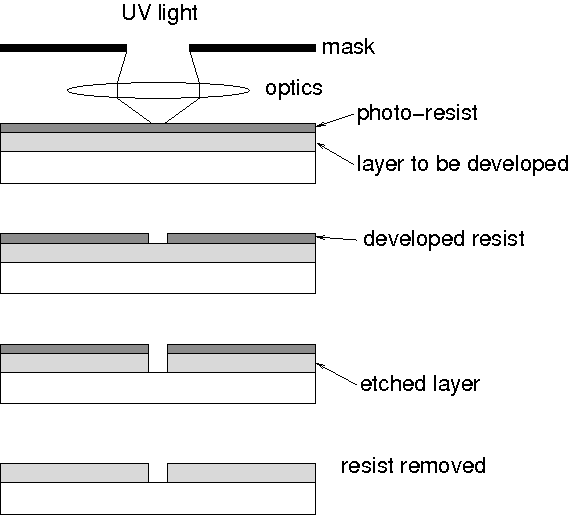}}
\caption{The photolithographic process.  Ultraviolet light, diffracted
  by a mask, forms an image on the photoresist.  The exposed portion
  of the photoresist is removed, leaving openings.  Etching removes
  parts of the layer to be patterned.  After etching, the photoresist
  is removed.}
\label{process}
\end{figure}

The problem we address in this work is the inverse problem of
determining what mask is needed in order to remove a desired shape
in the photoresist.  The difficulty of producing a desired shape comes
from the fact that the UV light is diffracted at the mask.  Moreover,
the chemicals in the photoresist reacts nonlinearly to UV exposure --
only portions of the photoresist that have been exposed to a certain
level of intensity are removed in the bleaching process.

The nature of the present work is analytical.  Our goal is to formulate
mathematically well-posed problems for photolithography.  The methods
we use to prove well-posedness are constructive and may serve as a
foundation for a computational method.

Our investigation into photolithography is inspired by the work of
Cobb \cite{Cobb98} who was the first to approach this problem from the
point of view of optimal design which utilizes a physically-based
model.  This general approach was further developed by introducing
a level set method in \cite{Tuzel09}.  A different computational
approach which models the mask as a pixelated binary image can be found in
\cite{Poona07}.

The plan of the paper is as follows. In the first and preliminary
section, Section~\ref{model}, we develop the most basic model for
removal of the exposed photoresist.  We describe the inverse problem
to be solved.  This is followed by a discussion of the approximate problem
whose properties we intend to investigate in this work.
Section~\ref{prelsec} contains mathematical preliminaries needed
for our work.  We introduce the basic notation and recall
various results which will be useful for our analysis.
In particular, in Subsection~\ref{geomsubsec} we discuss the
geometry of masks or circuits and how to measure the distance between
two of them.  In Section~\ref{dirsec} we discuss the properties of the operator which
maps the mask into the circuit. Section~\ref{approxsec} provides an
analysis of the variational approach to the problem of the optimization of the mask and we prove a
convergence result for it, Theorem~\ref{finalteo}, in the framework of
$\Gamma$-convergence.

\section{Description of the inverse problem}\label{model}
This section is separated into three subsections.  First, we review some
basic facts about Fourier transforms and prove a result about approximation
of a Gaussian.  We follow this with a discussion of
the optics involved and a model for photolithography.  In the final subsection
we describe the inverse problem and its approximation.

\subsection{Fourier transform and approximation of Gaussians}
We first set some notation and describe a few preliminary results.
For every $x\in\mathbb{R}^2$, we shall set $x=(x_1,x_2)$, where $x_1$ and $x_2\in\mathbb{R}$. For every $x\in\mathbb{R}^2$ and $r>0$, we shall denote by $B_r(x)$ the open ball
in $\mathbb{R}^2$ centered at $x$ of radius $r$.
Usually we shall write $B_r$ instead of $B_r(0)$.
We recall that, for any set $E\subset \mathbb{R}^2$, we denote by $\chi_E$
its characteristic function, and for any $r>0$,
$B_r(E)=\bigcup_{x\in E}B_r(x)$.

For any $f\in\mathcal{S}'(\mathbb{R}^2)$, the space of tempered distributions,
we denote by $\hat{f}$ its \emph{Fourier transform}, which, if
$f\in L^1(\mathbb{R}^2)$, may be written as
$$\hat{f}(\xi)=\int_{\mathbb{R}^2}f(x)\rme^{-\rmi \xi\cdot x}\rmd x,\quad \xi\in\mathbb{R}^2.$$
We recall that $f(x)=(2\pi)^{-2}\hat{\hat{f}}(-x)$, that is, when also $\hat{f}\in L^1(\mathbb{R}^2)$,
$$f(x)=\frac{1}{(2\pi)^2}\int_{\mathbb{R}^2}\hat{f}(\xi)\rme^{\rmi \xi\cdot x}\rmd \xi,\quad x\in\mathbb{R}^2.$$

If $f$ is a radial function, that is $f(x)=\phi(|x|)$ for any $x\in \mathbb{R}^2$, then
$$\hat{f}(\xi)=2\pi \mathcal{H}_0(\phi)(|\xi|),\quad \xi\in\mathbb{R}^2,$$
where
$$\mathcal{H}_0(\phi)(s)=\int_0^{+\infty}rJ_0(sr)\phi(r)\rmd r,\quad s\geq 0,$$
is the Hankel transform of order $0$, $J_0$ being the Bessel function of order $0$, see for instance \cite{Col}.

We denote the Gaussian distribution by
$G(x)=(2\pi)^{-1}\rme ^{-|x|^2/2}$, $x\in \mathbb{R}^2$, and let us note that
$\hat{G}(\xi)=\rme ^{-|\xi|^2/2}$, $\xi\in\mathbb{R}^2$. Moreover, $\|G\|_{L^1(\mathbb{R}^2)}=1$. Furthermore if $\delta_0$ denotes the Dirac delta centered at $0$, we have $\widehat{\delta_0}\equiv 1$, therefore $(2\pi)^{-2}\hat{1}=\delta_0$.

For any function $f$ defined on $\mathbb{R}^2$ and any positive constant $s$, we denote $f_s(x)=s^{-2}f(x/s)$,
$x\in \mathbb{R}^2$.
We note that $\|f_s\|_{L^1(\mathbb{R}^2)}=\|f\|_{L^1(\mathbb{R}^2)}$ and $\widehat{f_s}(\xi)=\hat{f}(s\xi)$, $\xi\in\mathbb{R}^2$.

We conclude these preliminaries with the following integrability result for the Fourier transform and its applications.

\begin{teo}\label{intFourierteo}
There exists an absolute constant $C$ such that the following estimate holds
$$\|\hat{f}\|_{L^1(\mathbb{R}^2)}\leq C\|f\|_{W^{2,1}(\mathbb{R}^2)}.$$
\end{teo}

\proof{.} This result is contained in Theorem A in \cite{Kol} and it is based on previous analysis done in \cite{Pel-Woj}.\cvd

\bigskip

We recall that a more detailed analysis on conditions for which integrability of the Fourier transform holds may be found in \cite{Taib}. However the previous result is simple to use and it is enough for our purposes, in particular for proving the following lemma.

\begin{lem}\label{approxlemma}
For any $\tilde{\delta}>0$ there exist a constant $s_0$, $0<s_0\leq 1$, and a radial function $\hat{T}\in C_0^{\infty}(\mathbb{R}^2)$ such that $\hat{T}\equiv 1$ on $B_{s_0}$ and, if we call
$T=(2\pi)^{-2}\hat{\hat{T}}$, then $T\in W^{2,1}(\mathbb{R}^2)$ and
$$\|T-G\|_{W^{1,1}(\mathbb{R}^2)}\leq \tilde{\delta}.$$
\end{lem}

\proof{.}
We sketch the proof of this result. Let us consider the following cut-off function
$\phi\in C^{\infty}(\mathbb{R})$ such that $\phi$ is nonincreasing,  $\phi\equiv 1$ on$(-\infty,0]$ and $\phi\equiv 0$ on $[1,+\infty)$.

We define a function $\hat{\tilde{T}}$ as follows
$$\hat{\tilde{T}}(x)=\phi(|x|-1)+(1-\phi(|x|-1))\widehat{G_{s_0}}(x)\phi(|x|-b),\quad x\in\mathbb{R}^2,$$
for suitable constants $s_0$, $0<s_0\leq 1$ and $b\geq 2$.
We call $\tilde{T}=(2\pi)^{-2}\hat{\hat{\tilde{T}}}$.

Then lengthy but straightforward computations, with the aid of Theorem~\ref{intFourierteo}, allow us to prove that for some $s_0$ small enough and for some $b=s_0^{-1}b_0$, with $b_0$ large enough, we have
$$\|\tilde{T} -G_{s_0}\|_{W^{1,1}(\mathbb{R}^2)}\leq \tilde{\delta}.$$

Then, let $\hat{T}(x)=\hat{\tilde{T}}(x/s_0)$, $x\in\mathbb{R}^2$, so that
$T=\tilde{T}_{1/s_0}$, or equivalently $\tilde{T}=T_{s_0}$.
Therefore
$$\|T_{s_0}-G_{s_0}\|_{W^{1,1}(\mathbb{R}^2)}\leq \tilde{\delta}.$$
By a simple rescaling argument we have that $\hat{T}$ satisfies the required properties.
Furthermore, by this construction, we may choose $\hat{T}$ such that it is radially nonincreasing, $\hat{T}\equiv 1$ on $B_{s_0}$ and it decays to zero in a suitable smooth, exponential way.\cvd

\subsection{A model of image formation}

We are now in the position to describe the model we shall use.
The current industry standard for modeling the optics is based on
Kirchhoff approximation.  Under this approximation,
the light source at the mask is on where the mask is open, and off
otherwise (see Figure 1).  Propagation through the lenses can be
calculated using Fourier optics.  It is further assumed that the image
plane, in this case the plane of the photoresist, is at the focal
distance of the optical system.  If there were no
diffraction, a perfect image of the mask would be formed on the image
plane. Diffraction, together with partial coherence of the light
source, acts to distort the formed image.

The mask, which we mention consists of cut-outs, is represented as a 
binary function, i.e., it is a
characteristic function of the cut-outs.  Suppose that $D$ represents
the cut-outs, then the mask is given by
\[
m(x) = \chi_D(x).
\]
The image is the light intensity on the image plane. This is given by
\cite{Pati97}
\begin{equation}\label{Hopkins}
I(x) = \int_{\mathbb{R}^2} \int_{\mathbb{R}^2} m(\xi)K(x-\xi) J(\xi-\eta) K(x-\eta) m(\eta) \rmd\xi \rmd\eta,\quad x\in\mathbb{R}^2. 
\end{equation}
In the above expression the kernel $K(\cdot)$ is called the \emph{coherent point spread function}
and describes the optical system.  For an optical system with a
circular aperture, once the wavenumber of the light
used, $k>0$, has
been chosen, the kernel depends on a single parameter called the Numerical
Aperture, $\text{NA}$. Notice that the wavelength
is $\lambda=2\pi/k$.
Let us recall that the so-called Jinc function is defined as
\[
\textrm{Jinc}(x)=\frac{J_1(|x|)}{2\pi|x|},\quad x\in\mathbb{R}^2,
\]
where $J_1$ is the Bessel function of order 1. We notice that in the 
Fourier space, see for instance \cite[page~14]{Erd},
$$\widehat{\textrm{Jinc}}(\xi)=\chi_{B_1}(\xi),\quad \xi\in\mathbb{R}^2.$$

If we denote by $s=(k\text{NA})^{-1}$, then
the kernel
is usually modeled as follows
\[
K(x)=\textrm{Jinc}_s(x)=\frac{k\text{NA}}{2\pi}\frac{J_1(k \text{NA}|x|)}{|x|},\quad x\in\mathbb{R}^2,
\]
therefore
$$\hat{K}(\xi)=\chi_{B_1}(s\xi)=\chi_{B_{1/s}}(s\xi)=\chi_{B_{k\text{NA}}}(\xi),\quad \xi\in\mathbb{R}^2.$$
If $\text{NA}$ goes to $+\infty$, that is $s\to 0^+$, then $\hat{K}$ converges pointwise to $1$, thus $K$ approximates in a suitable sense the Dirac delta.

For technical reasons, we shall consider a slightly different coherent point spread function $K$. Let us fix a positive constant $\tilde{\delta}$, to be chosen later.
We shall replace the characteristic function $\chi_{B_1}$, the Fourier transform of the Jinc function, with the function $\hat{T}(s_0\xi)$, $\xi\in\mathbb{R}^2$, with $\hat{T}$ and $s_0$ as in Lemma~\ref{approxlemma}.
Therefore $\hat{T}(s_0\cdot)$ is a radial function that it is still identically equal to $1$ on $B_1$, it is still compactly supported, it is nonincreasing with respect to the radial variable and it decays to zero in a smooth, exponential way.
Its Fourier transform is $T_{s_0}$ and we shall assume that
\begin{equation}\label{Kdef}
K(x)=(T_{s_0})_s(x)=T_{ss_0}(x),\quad x\in\mathbb{R}^2,
\end{equation}
where again $s=(k\text{NA})^{-1}$. Also in this model, if $\text{NA}$ goes to $+\infty$, that is $s\to 0^+$, then $\hat{K}$ converges pointwise to $1$, thus $K$ approximates in a suitable sense the Dirac delta.

The function $J(\cdot)$ is called the \emph{mutual
intensity function}.  If the illumination is fully coherent, $J\equiv 1$.
In practice, illumination is never fully coherent and is parametrized
by a \emph{coherency coefficient} $\sigma$.  A typical model for $J$ is
\begin{equation}\label{Jdef}
J(x)=2\frac{J_1(k \sigma \text{NA}|x|)}{k\sigma\text{NA}|x|}=\pi \textrm{Jinc}(k\sigma\text{NA}|x|),\quad x\in\mathbb{R}^2.
\end{equation}
Thus,
$$\frac{1}{(2\pi)^2}\hat{J}(\xi)=\frac{1}{\pi(k\sigma\text{NA})^2}\chi_{B_{k\sigma\text{NA}}}(\xi),\quad\xi\in\mathbb{R}^2,$$
that, as $\sigma\to 0^+$, converges, in a suitable sense, to the Dirac delta. Therefore
full coherence is achieved for $\sigma\to 0^+$. In fact, if $\sigma\to 0^+$, $J$ converges to $1$ uniformly on any compact subset of $\mathbb{R}^2$. 
The equation (\ref{Hopkins}) is often referred to as the Hopkins areal
intensity representation. As it will become apparent from the analysis developed in the paper, the value of $s$ is related to the scale of details that the manufacturing of the mask allows, thus in turn to the scale of details of the desired circuit. Therefore,
we typically consider $k\text{NA}\gg 1$, that is $s\ll 1$, and
$k\sigma\text{NA}\ll 1$.

\subsection{The inverse problem and its approximation}
The photoresist material responds to the intensity of the image.  When
intensity at the photoresist goes over a certain threshold, it is
then considered exposed and can be removed.  Therefore, the
exposed pattern, given a mask $m(x)$, is
\begin{equation}\label{exposed}
\Omega = \{ x\in\mathbb{R}^2\,:\ I(x) > h \},
\end{equation}
where $h$ is the exposure threshold.  Clearly, $\Omega$ depends on the
mask function $m(x)$, which we recall is given by the characteristic
function of $D$ representing the cut-outs, that is $\Omega=\Omega(D)$.  In photolithography, we have a
desired exposed pattern which we wish to achieve.  The inverse problem
is to find a mask that achieves this desired exposed pattern.
Mathematically, this cannot, in general, be done.  Therefore, the
inverse problem must be posed as an optimal design problem.

Suppose the desired pattern is given by $\Omega_0$.  We pose the
minimization problem
\begin{equation}\label{design_prob}
\displaystyle{\min_{D\in{\mathcal A}}} \; d(\Omega(D), \Omega_0).
\end{equation}
The distance function $d(\cdot,\cdot)$ will be discussed in detail
below.  The admissible set $\mathcal A$ is our search space, and needs to
be defined carefully as well. 

Instead of solving (\ref{design_prob}), we pose a variational problem for
a function $u$ (instead of the mask $D$).  We will show below that this problem is well-posed and that
as the approximation parameter is set to zero, we recover the solution of (\ref{design_prob}) under
a perimeter penalization.

Instead of dealing with the characteristic function $\chi_D(x)$ which represents
the mask, we will work with a phase-field function $u$ which takes on values of 0 and 1 with
smooth transitions.  Thus, the intensity in (\ref{exposed}) is calculated with $u$ instead of $m=\chi_D$ in \eqref{Hopkins},
so $I$ is a function of $u$.  At this point, we will not be precise about the space of
functions to which $u$ belongs.  To force $u$ to take on values of mostly 0 and 1, we
introduce the Mordica-Mortola energy
\[
P_\varepsilon(u) = 
\displaystyle{\frac{1}{\varepsilon}\int W(u)+\varepsilon\int|\nabla u|^2},
\]
where $\displaystyle{W(t)=9t^2(t-1)^2}$ is a double-well potential.
We will regularize the problem of minimizing the distance between the target pattern
and the exposed region by this energy.

Then we relax the hard threshold in defining the exposed 
region $\Omega$ in (\ref{exposed}).
Let $\phi(t)$ be a $C^{\infty}$ nondecreasing approximate Heaviside function 
with values $\phi(t \leq -1/2)=0$ and $\phi(t\geq 1/2)=1$.  The function
$$\Phi_{\eta}(u)=\phi\left(\frac{I(u)-h}{\eta} \right) $$
will be 1 where the intensity $I \geq h+\eta/2$.  A sigmoidal threshold function is
employed in the computational work in \cite{Poona07}.

Now we consider the distance function between $\Omega$ and $\Omega_0$ in (\ref{design_prob}).  Let
\begin{equation}
d=d(\Omega,\Omega_0)=\int|\chi_{\Omega}-\chi_{\Omega_0}|+\left|P(\Omega)-P(\Omega_0)\right|,
\end{equation}
where $\chi_\Omega$ is the characteristic function of the set $\Omega$ and
$P(\Omega)$ is the perimeter of the region $\Omega$.  To approximate this distance
function, we replace it by
\[
d_\eta(u,\Omega_0) = \int | \Phi_\eta(u) - \chi_{\Omega_0} | + 
\left| \int | \nabla (\Phi_\eta(u))| - P(\Omega_0) \right|.
\]
The characteristic function of $\Omega$ is replaced by the smooth threshold function while
its perimeter is replaced by the TV-norm of the function.  

The approximate problem we shall solve is 
\[
F_{\varepsilon}(u)= d_{\eta(\varepsilon)}(u,\Omega_0) + b P_\varepsilon(u) \rightarrow \min .
\]
The remainder of the paper is an analytical study of this minimization problem.  We will show
that it is well-posed, and that in the limit $\varepsilon \rightarrow 0^+$, we recover the
solution of the original problem (\ref{design_prob}) under a perimeter penalization.

\section{Mathematical preliminaries}\label{prelsec}
\newcommand{\mcd}{\mathcal{D}}
By $\mathcal{H}^1$ we denote the $1$-dimensional Hausdorff measure and by
$\mathcal{L}^2$ we denote the $2$-dimensio\-nal Lebesgue measure. We recall that,
if $\gamma\subset\mathbb{R}^2$ is a smooth curve,
then $\mathcal{H}^1$ restricted to
$\gamma$ coincides with its arclength.
For any Borel $E\subset\mathbb{R}^2$ we denote $|E|=\mathcal{L}^2(E)$.

Let $\mcd$ be a bounded open set contained in $\mathbb{R}^2$, with boundary $\partial \mcd$.
We say that $\mcd$ has a \emph{Lipschitz boundary}
if for every $x=(x_1,x_2)\in\partial\mcd$ there exist a Lipschitz
function $\varphi:\mathbb{R}\to\mathbb{R}$ and a positive constant $r$
such that for any $y\in B_r(x)$ we have, up to a rigid transformation,
$$
y=(y_1,y_2)\in\mcd\quad \text{if and only if}\quad  y_2<\varphi(y_1).
$$
We note that $\mcd$ has a finite number of connected components, whereas
$\partial \mcd$ is formed by a finite number of rectifiable Jordan curves, therefore
$\mathcal{H}^1(\partial\mcd)=\mathrm{length}(\partial\mcd)<+\infty$.

We recall some basic notation and properties
of functions of bounded variation and sets of finite perimeter. For a more comprehensive treatment of
these subjects see, for instance, \cite{Amb e Fus e Pal, Eva e Gar, Giu}.

Given a bounded open set $\mcd\subset \mathbb{R}^2$,
we denote by $BV(\mcd)$ the Banach space of \emph{functions of bounded
variation}. We recall that $u\in BV(\mcd)$ if and only if  
$u\in L^1(\mcd)$ and its distributional derivative $Du$ is a bounded
vector
measure. We endow $BV(\mcd)$ with the standard norm as follows. Given
$u\in BV(\mcd)$, we denote by $|Du|$ the total variation of its
distributional derivative and
we set $\|u\|_{BV(\mcd)}=\|u\|_{L^1(\mcd)}+|Du|(\mcd)$. We shall call
$P(u,\mcd)=|Du|(\mcd)$.
We recall that whenever $u\in W^{1,1}(\mcd)$,
then $u\in BV(\mcd)$ and $|Du|(\mcd)=\int_{\mcd}|\nabla u|$, therefore
$\|u\|_{BV(\mcd)}=\|u\|_{L^1(\mcd)}+\|\nabla u\|_{L^1(\mcd)}=\|u\|_{W^{1,1}(\mcd)}$.

We say that a sequence of $BV(\mcd)$ functions $\{u_h\}_{h=1}^{\infty}$
\emph{weakly}$^*$ \emph{converges} in $BV(\mcd)$ to $u\in BV(\mcd)$ if and only if
$u_h$ converges to $u$ in $L^1(\mcd)$ and $Du_h$
weakly$^*$ converges to $Du$ in $\mcd$, that is
\begin{equation}\label{weakstarconv}
\lim_{h}\int_{\mcd}v \rmd Du_h=
\int_{\mcd}v \rmd Du\quad\text{for any }v\in C_0(\mcd).
\end{equation}
By Proposition~3.13 in \cite{Amb e Fus e Pal}, we have that if
a sequence of $BV(\mcd)$ functions $\{u_h\}_{h=1}^{\infty}$ is bounded in $BV(\mcd)$ and converges to $u$ in $L^1(\mcd)$, then
$u\in BV(\mcd)$ and $u_h$ converges to $u$
weakly$^*$ in $BV(\mcd)$.

We say that a sequence of $BV(\mcd)$ functions $\{u_h\}_{h=1}^{\infty}$
\emph{strictly converges} in $BV(\mcd)$ to $u\in BV(\mcd)$ if and only if
$u_h$ converges to $u$ in $L^1(\mcd)$ and $|Du_h|(\mcd)$
converges to $|Du|(\mcd)$. Indeed,
$$d_{st}(u,v)=\int_{\mcd}|u-v|+\big||Du|(\mcd)-|Dv|(\mcd) \big|$$
is a distance on $BV(\mcd)$ inducing the strict convergence.
We also note that strict convergence implies weak$^*$ convergence.

Let $\mcd$ be a bounded open set with Lipschitz boundary.
A sequence of $BV(\mcd)$ functions $\{u_h\}_{h=1}^{\infty}$
such that $\sup_h\|u_h\|_{BV(\mcd)}<+\infty$ admits a subsequence
converging weakly$^*$ in $BV(\mcd)$ to a function $u\in BV(\mcd)$, see for instance Theorem~3.23 in \cite{Amb e Fus e Pal}.
As a corollary, we infer that for any $C>0$ the set
$\{u\in BV(\mcd)\,:\ \|u\|_{BV(\mcd)}\leq C\}$
is a compact subset of $L^1(\mcd)$.

For any fixed constant $R>0$, with a slight abuse of notation, we shall identify
$L^1(B_R)$ with the set $\{u\in L^1(\mathbb{R}^2)\,:\ u=0\text{ a.e. outside }B_R\}$. 

Let $E$ be a bounded 
Borel set contained in $B_R\subset \mathbb{R}^2$. We shall denote by $\chi_E$ its characteristic function. We notice that 
$E$ is compactly
contained in $B_{R+1}$, which we shall denote by $E\Subset B_{R+1}$.
We say that $E$ is a \emph{set of finite perimeter} if
$\chi_E$ belongs to $BV(B_{R+1})$ and we call the number
$P(E)=|D\chi_E|(B_{R+1})$ its \emph{perimeter}.
Analogously, for any $u\in L^1(B_R)\cap BV(B_{R+1})$, we shall denote
$P(u,B_{R+1})=|Du|(B_{R+1})$. Obviously, if
$u=\chi_E$, then $P(u,B_{R+1})=P(E)$.

Let us further remark that the intersection of two sets of finite perimeter is
still a set of finite perimeter. Moreover,
whenever $E$ is open and $\mathcal{H}^1(\partial E)$ is finite, then $E$ is a set of finite
perimeter, see for instance \cite[Section~5.11, Theorem~1]{Eva e Gar}.
Therefore a bounded open set
$\mcd$ with Lipschitz boundary
is a set of finite perimeter and its perimeter $P(\mcd)$ coincides with
$\mathcal{H}^1(\partial\mcd)$.

\subsection{$\Gamma$-convergence approximation of the perimeter functional}\label{ModMorsubsec}

Let us introduce the following, slightly different, version of a $\Gamma$-convergence result due to Modica and Mortola, \cite{Mod e Mor}. We shall follow the notation and proofs contained in
\cite{Bra}. We begin by setting some notation. For the definition and properties of $\Gamma$-convergence we refer to \cite{DaM}.

For any bounded open set $\mcd\subset\mathbb{R}^2$, with a slight abuse of notation,
we identify $W^{1,p}_0(\mcd)$, $1<p<+\infty$, with the subset of $W^{1,p}(\mathbb{R}^2)$ functions $u$ such that $u$ restricted
to $\mcd$ belongs to $W^{1,p}_0(\mcd)$ and $u$ is equal to $0$ almost everywhere outside $\mcd$.
Let us assume that for some positive constant $R$ we have $\mcd\subset B_R$.
We recall that any function in $L^1(\mcd)$ is extended to zero outside $\mcd$
and the same procedure is used for $L^1(B_R)$.
Therefore, with this slight abuse of notation, $L^1(\mcd)\subset L^1(B_R)$.
Throughout the paper, for any $p$, $1\leq p\leq +\infty$, we shall denote its conjugate exponent by $p'$, that is $p^{-1}+(p')^{-1}=1$.

\begin{teo}\label{Mod-Morteo}
Let $\mcd\subset B_R\subset\mathbb{R}^2$ be a bounded open set with Lipschitz boundary. Let us also assume that
$\mcd$ is convex.

Let $1<p<+\infty$ and $W:\mathbb{R}\to[0,+\infty)$ be a continuous function such that
$W(t)=0$ if and only if $t\in\{0,1\}$. Let $c_p=(\int_0^1(W(s))^{1/p'}\mathrm{d}s)^{-1}$.

For any $\varepsilon>0$ we define the functional
$P_{\varepsilon}:L^1(\mathbb{R}^2)\to [0,+\infty]$ as follows
\begin{equation}\label{modmordef}
P_{\varepsilon}(u)=\left\{
\begin{array}{ll}
\displaystyle{\frac{c_p}{p'\varepsilon}\int_{\mcd}W(u)+\frac{c_p\varepsilon^{p-1}}{p}\int_{\mcd}|\nabla u|^p}&\text{if }u\in W^{1,p}_0(\mcd),\\
\vphantom{\displaystyle{\int}}+\infty&\text{otherwise}.
\end{array}
\right.
\end{equation}

Let $P:L^1(\mathbb{R}^2)\to [0,+\infty]$ be such that
\begin{equation}\label{Pdef}
P(u)=
\left\{
\begin{array}{ll}
\vphantom{\displaystyle{\int}}P(u,B_{R+1})
&\text{if }u\in BV(B_{R+1}),\ u\in\{0,1\}\text{ a.e.},\\ &\quad\text{and }
u=0\text{ a.e. outside }\mcd,\\
\vphantom{\displaystyle{\int}}+\infty&\text{otherwise}.
\end{array}
\right.
\end{equation}

Then $P=\Gamma\textrm{-}\!\lim_{\varepsilon\to 0^+}P_{\varepsilon}$ with respect to the $L^1(\mathbb{R}^2)$ norm.
\end{teo}

\begin{oss}
We observe that $P(u)=P(E)$ if $u=\chi_E$ where $E$ is a set of finite perimeter contained in $\overline{\mcd}$ and $P(u)=+\infty$ otherwise.

Furthermore, we note that the result does not change if in the definition of $P_{\varepsilon}$ we set $P_{\varepsilon}(u)=+\infty$ whenever $u$ does not satisfy the
constraint
\begin{equation}
0\leq u\leq 1\text{ a.e. in }\mcd.
\end{equation}
\end{oss}

\proof{.} We sketch the proof following that of Theorem~4.13 in \cite{Bra}. In fact,
the only difference with respect to that theorem is that we assume $\mcd$ convex and that we take $W^{1,p}_0(\mcd)$
instead of $W^{1,p}(\mcd)$ in the definition
of $P_{\varepsilon}$.

By Proposition~4.3 in \cite{Bra},  
we obtain that $P(u)\leq \Gamma\textrm{-}\!\liminf_{\varepsilon\to 0^+} P_{\varepsilon}(u)$ for any $u\in L^1(\mathbb{R}^2)$.
In order to obtain the $\Gamma\textrm{-}\!\limsup$ inequality, we follow the procedure described in
Section~4.2 of \cite{Bra}.
It would be enough to construct
$\mathcal{M}\subset L^1(\mathbb{R}^2)$ such that the following two conditions are satisfied.
First, we require that, for any $u\in L^1(\mathbb{R}^2)$ such that $P(u)<+\infty$, there exists a sequence
$\{u_j\}_{j=1}^{\infty}$ such that
$u_j\in \mathcal{M}$, for any $j\in\mathbb{N}$, $u_j\to u$ in $L^1(\mathbb{R}^2)$ as $j\to\infty$, and $P(u)=\lim_jP(u_j)$.
Second, for any $u\in\mathcal{M}$, $\Gamma\textrm{-}\!\limsup_{\varepsilon\to 0^+} P_{\varepsilon}(u)\leq P(u)$.

We choose $\mathcal{M}=\{u=\chi_E\,:\ E\Subset\mcd,\ E\text{ of class }C^{\infty}\}$.
The second property follows by Proposition~4.10 in \cite{Bra}. As far as the first property is concerned, this can be
obtained by following the proof of Theorem~1.24 in \cite{Giu}. That theorem states that any bounded set of finite perimeter
$E$ can be approximated by a sequence of $C^{\infty}$ sets $\{E_j\}_{j=1}^{\infty}$ such that, as $j\to \infty$,
$\int_{\mathbb{R}^2}|\chi_{E_j}-\chi_E|\to 0$ and $P(E_j)\to P(E)$. If we assume that $E\subset \overline{\mcd}$,
and that $\mcd$ is convex,  by choosing in the proof of Theorem~1.24 in \cite{Giu} a value of $t$ 
satisfying
$1/2<t<1$, we obtain that the sets $E_j$ are also compactly contained in $\mcd$, for any $j\in\mathbb{N}$.\cvd

\bigskip

Also the following result, due to Modica, \cite{Mod}, will be useful.

\begin{prop}\label{modcompprop}
For any $C>0$, let us take $1<p<+\infty$ and any $\varepsilon>0$, and let us define
$$A_C=\{u\in L^1(\mathbb{R}^2)\,:\ 0\leq u\leq 1\text{ a.e. and }P_{\varepsilon}(u)\leq C\}.$$
Then $A_C$ is precompact in $L^1(\mathbb{R}^2)$.
\end{prop}

\proof{.} We repeat, for the reader's convenience, the arguments developed in \cite{Mod}.
Clearly $A_C$ is a bounded subset of $L^1(\mcd)$.
Let $\{u_n\}_{n=1}^{\infty}$ be a sequence in $A_C$. We need to prove that there exists a subsequence converging in $L^1(\mcd)$.
For any $t$, $0\leq t\leq 1$, let
$\phi(t)=\int_0^t(W(s))^{1/p'}\mathrm{d}s$. For any $n\in\mathbb{N}$,
we define $v_n=\phi(u_n)$ and we observe that
$0\leq v_n\leq \phi(1)$ almost everywhere. Therefore,
the functions $v_n$, $n\in\mathbb{N}$, are 
uniformly bounded in $L^{\infty}(\mcd)$ and, consequently, in $L^1(\mcd)$. Furthermore, since $\phi$ is a $C^1$ function,
with bounded $C^1$ norm, then $Dv_n=\phi'(u_n)Du_n=W^{1/p'}(u_n)Du_n$.
Therefore,
$$\int_{\mcd}|Dv_n|=\int_{\mcd}|W^{1/p'}(u_n)||Du_n|\leq P_{\varepsilon}(u_n)/c_p.$$
We infer that there exists a subsequence $\{v_{n_k}\}_{k=1}^{\infty}$ converging, as $k\to\infty$, to a function $v_0$ in $L^1(\mcd)$
and almost everywhere.
Let $\psi$ be the inverse function of $\phi$ and let $u_0=\psi(v_0)$. We observe that
$\psi$ is bounded and uniformly continuous on $[0,\phi(1)]$, hence we conclude that,
as $k\to\infty$, $u_{n_k}$ converges to $u_0$ in $L^1(\mcd)$.\cvd

\begin{oss}\label{compactnessoss}
With the same proof, we can show the following. Let us consider
any family $\{u_{\varepsilon}\}_{0<\varepsilon\leq \varepsilon_0}$ such that,
for some positive constant $C$ and
for any $\varepsilon$, $0<\varepsilon\leq\varepsilon_0$,
we have $0\leq u_{\varepsilon}\leq 1$ almost everywhere and
$P_{\varepsilon}(u_{\varepsilon})\leq C$.
Then $\{u_{\varepsilon}\}_{0<\varepsilon\leq \varepsilon_0}$
is precompact in $L^1(\mathbb{R}^2)$.
\end{oss}


\subsection{Convolutions}

We recall that, for any two functions $f$ and $g$ defined on $\mathbb{R}^2$,
we define the \emph{convolution} of $f$ and $g$, $f\ast g$, as follows
$$(f\ast g)(x)=\int_{\mathbb{R}^2}f(x-y)g(y)\rmd y=\int_{\mathbb{R}^2}f(y)g(x-y)\rmd y,\quad x\in\mathbb{R}^2,$$
whenever this is well-defined.

The following classical properties of convolutions will be used. First convolution is commutative. Second, as a consequence of Young inequality we have
the following result about integrability and regularity of convolutions.

\begin{prop}\label{convprop}
Let $1\leq r,\ p,\ q\leq +\infty$ be such that $1+\frac{1}{r}=\frac{1}{p}+\frac{1}{q}$, and let $n=0,1,2,\ldots$.

Let $1\leq q<+\infty$,
let $f\in L^{q}(\mathbb{R}^2)$ and let $g\in W^{n,p}(\mathbb{R}^2)$. Then
$h=f\ast g\in W^{n,r}(\mathbb{R}^2)$ and there exists a constant $C$, depending on $n$, $p$, $q$ and $r$ only, such that
$$\|h\|_{W^{n,r}(\mathbb{R}^2)}\leq C\|f\|_{L^q(\mathbb{R}^2)}\|g\|_{W^{n,p}(\mathbb{R}^2)}.$$

Let $q=+\infty$ and let  $f\in L^{\infty}(\mathbb{R}^2)$, with compact support. If 
$g\in W^{n,1}(\mathbb{R}^2)$, then
$h=f\ast g\in W^{n,\infty}(\mathbb{R}^2)$ and there exists a constant $C$, depending on $n$ only, such that
$$\|h\|_{W^{n,\infty}(\mathbb{R}^2)}\leq C\|f\|_{L^{\infty}(\mathbb{R}^2)}\|g\|_{W^{n,1}(\mathbb{R}^2)}.$$


If $f\in L^1(\mathbb{R}^2)$ and $g\in L^{\infty}(\mathbb{R}^2)$, then
$h=f\ast g\in L^{\infty}(\mathbb{R}^2)$ and it holds
$\|h\|_{L^{\infty}(\mathbb{R}^2)}\leq \|f\|_{L^{\infty}(\mathbb{R}^2)}\|g\|_{L^1(\mathbb{R}^2)}$. Furthermore, if $g$ is uniformly continuous and 
$\omega_g$ denotes its modulus of continuity, then $h$ is also uniformly continuous and 
$$\omega_h\leq \|f\|_{L^1(\mathbb{R}^2)}\omega_g.$$

Finally, let $f\in L^1(\mathbb{R}^2)$ and let $g\in C^{n,\alpha}(\mathbb{R}^2)$,
for some $\alpha$, $0<\alpha\leq 1$. Then $h\in C^{n,\alpha}(\mathbb{R}^2)$ 
and there exists a constant $C$, depending on $n$ and $\alpha$ only, such that
$$\|h\|_{C^{n,\alpha}(\mathbb{R}^2)}\leq C\|f\|_{L^1(\mathbb{R}^2)}\|g\|_{C^{n,\alpha}(\mathbb{R}^2)}.$$
\end{prop}

\subsection{The geometry of masks and circuits}\label{geomsubsec}

In this subsection we investigate the following two questions, namely
what are reasonable assumptions on the geometry of the mask $D$ and how to measure the distance between
the constructed circuit $\Omega$ and the desired one $\Omega_0$. We begin with the following definition. During this subsection, in most cases proofs will be omitted and left to the reader.

For given positive constants $r$ and $L$, we say that a bounded open set
$\Omega\subset \mathbb{R}^2$
is \emph{Lipschitz} or $C^{0,1}$ \emph{with constants} $r$ \emph{and} $L$ if 
for every $x\in\partial\Omega$ there exists a Lipschitz
function $\varphi:\mathbb{R}\to\mathbb{R}$, with Lipschitz constant bounded by $L$, such that for any $y\in B_r(x)$, and up to a rigid transformation,
\begin{equation}\label{Lipdomain}
y=(y_1,y_2)\in\Omega\quad \text{if and only if}\quad  y_2<\varphi(y_1).
\end{equation}
Without loss of generality, we may always assume that $x=(0,0)$ and $\varphi(0)=0$.
We shall always denote by $e_1$ and $e_2$ the vectors of the canonical bases.
Clearly the orientation of the canonical bases may
vary depending on $x\in\partial\Omega$.

We shall also use the following notation. 
There exist positive constants $\delta_1\leq 1/2$, $\delta_2\leq \delta_1$ and $m_1\leq 1$, all of them depending on $L$ only, such that the following holds.
For any $x\in \partial\Omega$ and for any $\delta>0$, let $M_{\delta}(x)=
\{y\,:\ |y_1|\leq\delta r,\ y_2=\varphi(y_1)\}$ and
$N_{\delta}(x)=\{y\,:\ |y_1|\leq\delta_1r,\ \varphi(y_1)-\delta r\leq y_2\leq\varphi(y_1)+\delta r\}$.
Then we assume that, for any $\delta$, $0<\delta\leq \delta_2$, the following properties hold. First,
$N_{\delta}(x)\subset B_{r/2}(x)$ (hence
$M_{\delta_1}(x)\subset B_{r/2}(x)$ as well).
Clearly $N_{\delta}(x)$ is contained in
$\overline{B}_{\delta r}(\partial\Omega)$, and we assume that $N_{\delta}(x)$
contains $\overline{B}_{m_1\delta r}(M_{\delta_1/2}(x))$ and that for any
$y\in \{y\,:\ |y_1|\leq\delta_1r/2,\ y_2=\varphi(y_1)\pm\delta r\}$, $y\not\in\overline{B}_{m_1\delta r}(\partial\Omega)$.

For any integer $k=1,2,\ldots$, any $\alpha$, $0< \alpha\leq 1$, and any positive constants $r$ and $L$,
we say that a bounded open set $\Omega\subset\mathbb{R}^2$
is $C^{k,\alpha}$ \emph{with constants} $r$ \emph{and} $L$ if 
for every $x\in\partial\Omega$ there exists a $C^{k,\alpha}$
function $\varphi:\mathbb{R}\to\mathbb{R}$, with $C^{k,\alpha}$ norm bounded by $L$, such that for any
$y\in B_r(x)$, and up to a rigid transformation, \eqref{Lipdomain} holds.
Without loss of generality, we may always assume that $x=(0,0)$ and $\varphi(0)=0$.

Let us fix three positive constants $r$, $L$ and $R$.
Let $\mathcal{A}^{0,1}(r,L,R)$ be the class of all bounded open sets, contained in $B_R\subset \mathbb{R}^2$,
which are Lipschitz with constants $r$ and $L$.
For any integer $k=1,2,\ldots$ and any $\alpha$, $0< \alpha\leq 1$, we denote with
$\mathcal{A}^{k,\alpha}(r,L,R)$ the class of all bounded open sets, contained in $B_R\subset \mathbb{R}^2$,
which are $C^{k,\alpha}$ with constants $r$ and $L$.

Since we shall identify open sets $D$ with their characteristic functions $\chi_{D}$, 
if $\mathcal{A}=\mathcal{A}^{0,1}(r,L,R)$,
(or $\mathcal{A}=\mathcal{A}^{k,\alpha}(r,L,R)$, respectively)
then, with a slight
abuse of notation,
$\mathcal{A}$ will also denote the subset of functions $u\in L^1(B_R)$ such that
$u=\chi_{D}$ for some $D\in \mathcal{A}$. Moreover,
we shall denote
$$A=\{u\in L^1(B_R)\,:\ 0\leq u\leq 1\text{ a.e. in }B_R\}$$
and, for any $\gamma>0$,
\begin{equation}
\mathcal{A}_{\gamma}=\{u\in A\,:\
\|u-\chi_{D}\|_{L^1(B_R)}\leq \gamma\text{ for some }D\in\mathcal{A}\}.
\end{equation}

Let us assume that
$\Omega_1$ and $\Omega_2$ belong to $\mathcal{A}^{0,1}(r,L,R)$.
There are several ways to define the distance between these two sets.
We shall describe four of them and study their relationships.
We let
\begin{align}
&d_1=d_1(\Omega_1,\Omega_2)=d_{H}(\overline{\Omega}_1,\overline{\Omega}_2);\\
&\tilde{d}_1=\tilde{d}_1(\Omega_1,\Omega_2)=d_{H}(\partial\Omega_1,\partial\Omega_2);\\
&d_2=d_2(\Omega_1,\Omega_2)=|\Omega_1\Delta\Omega_2|=\|\chi_{\Omega_1}-\chi_{\Omega_2}\|_{L^1(B_{R+1})};\\
&d_3=d_3(\Omega_1,\Omega_2)=d_2+|P(\Omega_1)-P(\Omega_2)|
=d_{st}(\chi_{\Omega_1},\chi_{\Omega_2}).\label{d3def}
\end{align}

Here $d_H$ denotes the Hausdorff distance, whereas we recall that
$P(\Omega)$ denotes the perimeter of $\Omega$ in $B_{R+1}$ and $d_{st}$ is the distance inducing strict convergence in $BV(B_{R+1})$.
First of all, we observe that all of these are distances. We now investigate their relationships.

We begin with the first two, $d_1$ and $\tilde{d}_1$, and we notice that
\begin{equation}
\text{if }d_1\leq r/4,\text{ then }d_1\leq \tilde{d}_1.
\end{equation}

There exists a constant $c$, $0<c\leq 1$, depending on $L$ only, such that
$$d_1\geq c \min\{r,\tilde{d}_1\}.$$
Therefore,
$$
\text{if }\tilde{d}_1\leq r,\text{ then }\tilde{d}_1\leq C d_1,
$$
where $C=1/c$. Furthermore, if $d_1\leq (c/2) r$, then $\tilde{d}_1$ must be less than or equal to $r$, so
$$
\text{if }d_1\leq (c/2)r,\text{ then }\tilde{d}_1\leq C d_1.
$$

Moreover, we can find a constant $c_1$, $0<c_1\leq 1$, depending on $L$ only,
such that
$$
\text{if }\tilde{d}_1\leq c_1r,\text{ then }d_1\leq \tilde{d}_1\leq C d_1.
$$
 
We conclude that we can find a constant $c_1$, $0<c_1\leq 1$,
 and a constant $C\geq 1$, both depending on $L$ only, such that
\begin{equation}\label{est1}
\text{if either }d_1\leq c_1r\text{ or }\tilde{d}_1\leq c_1r,\text{ then }d_1\leq \tilde{d}_1\leq C d_1.
\end{equation}

Since $d_1$ and $\tilde{d}_1$ are bounded by $2R$, we also have
\begin{equation}\label{est2}
\text{if both }d_1\geq c_1r\text{ and }\tilde{d}_1\geq c_1r,\text{ then }d_1\leq \frac{2R}{c_1r}\tilde{d}_1
\text{ and }\tilde{d}_1 \leq \frac{2R}{c_1r} d_1.
\end{equation}
We finally observe that the estimates \eqref{est1} and \eqref{est2} are essentially optimal.

Before comparing $d_1$ (or $\tilde{d}_1$) with $d_2$ and $d_3$, let us make the following remark on 
the lengths of $\partial\Omega_1$ and $\partial\Omega_2$.
If $\Omega$ is an open set which is Lipschitz with constants $r$ and $L$, then
for any integer $n\geq 0$, we have
\begin{equation}\label{length1}
\mathcal{H}^1(\partial\Omega\cap(\overline{B}_{(n+1)r}\backslash B_{nr}))\leq C(L)r(n+1).
\end{equation}
Here, a simple computation shows that we may choose
$C(L)=48 \sqrt{1+L^2}$.

Therefore, if we assume that $\Omega\subset B_R$ and $R\geq 10r$,
we may conclude that
\begin{equation}\label{length}
P(\Omega)\leq C_1(L)R^2/r,
\end{equation}
where $C_1(L)=\frac{1}{2}\left(\frac{11}{9}\right)^2C(L)$.

Moreover, there exist two constants $c_2$, $0< c_2\leq  c_1$, and $C_1>0$, depending on $L$ only, such that we have
\begin{equation}\label{neigh}
|\overline{B}_{d}(\partial\Omega)|\leq C_1\mathrm{length}(\partial\Omega)d\text{  for any }d\leq c_2r.
\end{equation}
Since
$$
\text{if }\tilde{d}_1\leq c_2r,\text{ then } d_2\leq \min\{|\overline{B}_{\tilde{d}_1}(\partial\Omega_2)|,|\overline{B}_{\tilde{d}_1}(\partial\Omega_2)| \},
$$
we obtain that
\begin{equation}\label{d2vsdtilde1}
\text{if }\tilde{d}_1\leq c_2r,\text{ then } d_2\leq C_1\min\{\mathrm{length}(\partial\Omega_1),\mathrm{length}(\partial\Omega_2)\}\tilde{d}_1.
\end{equation}

If $\tilde{d}_1\geq c_2r$, then
$d_2\leq \pi R^2\leq\frac{\pi R^2}{c_2r}\tilde{d}_1$. By \eqref{length},
we may conclude that
\begin{equation}\label{d2vsdtilde1bis}
d_2\leq \frac{C_2R^2}{r}\tilde{d}_1.
\end{equation}
Here $C_2$ depends on $L$ only. Moreover, up to changing the constants $c_2$, $C_1$ and $C_2$, \eqref{d2vsdtilde1} and
\eqref{d2vsdtilde1bis} still hold if we replace $\tilde{d}_1$ with $d_1$.

On the other hand, there exists a constant $c_3$, $0<c_3\leq \pi$, depending on $L$ only,
such that
$$d_2\geq c_3\min\{r^2,d_1^2\}.$$
We infer that either if $d_1\leq r$ or if $d_2\leq (c_3/2) r^2$, then
$d_1\leq C_3d_2^{1/2}$, where $C_3=1/c_3^{1/2}$. If $d_2\geq (c_3/2)r^2$,
then $d_1\leq 2R\leq \frac{4C_3^2}{r^2}Rd_2$
or, better,
$d_1\leq \frac{2\sqrt{2}C_3}{r}Rd_2^{1/2}$.
Summarizing, we have
\begin{equation}
\text{if }d_2\leq (c_3/2)r^2,\text{ then }d_1\leq C_3d_2^{1/2}
\end{equation}
and,
finally,
\begin{equation}
\text{if }d_2\geq (c_3/2)r^2,\text{ then }d_1\leq \frac{2\sqrt{2}C_3}{r}Rd_2^{1/2}.
\end{equation}
Clearly, up to suitably changing the constants $c_3$ and $C_3$, the last two estimates still hold if we
replace $d_1$ with $\tilde{d}_1$. We also remark that, as before, the estimates relating $d_1$, $\tilde{d}_1$ and $d_2$ are essentially optimal.

We have obtained that $d_1$, $\tilde{d}_1$ and $d_2$ are topologically equivalent distances.
About $d_2$ and $d_3$, obviously $d_2\leq d_3$, however the two distances are
not topologically equivalent. In fact we can find $\Omega$ and $\Omega_i$, $i\in\mathbb{N}$, open sets belonging to $\mathcal{A}^{0,1}(r,L,R)$, such that
$d_2(\Omega,\Omega_i)$ goes to zero as $i\to\infty$, whereas
$d_3(\Omega,\Omega_i)\geq c>0$ for any $i\in\mathbb{N}$. Therefore $d_3$ induces
a strictly finer topology than the one induced by $d_2$

An assumption that the mask is a bounded open set which is Lipschitz with
given constants $r$ and $L$ is reasonable from the manufacturing point of view as well as
from the mathematical point of view, by the 
following compactness result.

\begin{prop}\label{compactprop}
The set $\mathcal{A}^{0,1}(r,L,R)$
\textnormal{(}respectively $\mathcal{A}^{k,\alpha}(r,L,R)$, $k=1,2,\ldots$,
$0<\alpha\leq 1$\textnormal{)}
is compact with respect to the distance $d_1$.
\end{prop}

We remark that the same result holds with respect to the distances $\tilde{d}_1$ and
$d_2$. Furthermore, we obtain as a corollary that the set $\mathcal{A}_{\gamma}$
is closed with respect to the $L^1$ norm, for any $\gamma>0$.

The previous example shows that
compactness fails with respect to the distance $d_3$, at least for the Lipschitz case. On the other hand,
if $\Omega_1$ and $\Omega_2$ belong to $\mathcal{A}^{1,\alpha}(r,L,R)$,
with $0<\alpha<1$, then, following 
Lemma~2.1 in \cite{Ron99}, we can show that
\begin{equation}\label{d2vsd3}
|P(\Omega_1)-P(\Omega_2)|\leq C_4(\tilde{d}_1(\Omega_1,\Omega_2))^{\alpha/(2\alpha+2)},
\end{equation}
where $C_4$ depends on $r$, $L$, $R$ and $\alpha$ only.
We may conclude that in the $C^{k,\alpha}$ case, $k=1,2,\ldots$,
$0<\alpha\leq 1$, $d_3$ is topologically equivalent to the other three distances and that
Proposition~\ref{compactprop} holds also with respect to the distance $d_3$.

It is worthwhile to observe that, under some circumstances, the estimate
\eqref{d2vsd3} can be extended to the piecewise $C^{1,\alpha}$ case.
For example, typically we may assume that the desired circuit $\Omega_0$ belongs
to $\mathcal{A}^{0,1}(r,L,R)$.
Moreover, we assume that the boundary of $\Omega_0$ is composed by a finite number of closed segments
$I_i$, $i=1,\ldots,n$, which are pairwise internally disjoint and
whose lengths are greater than or equal to $2r$. Therefore, $\Omega_0$ is actually
a piecewise $C^{1,\alpha}$ open set. 
We shall show in Section~\ref{dirsec} that, under suitable assumptions on the mask $D$,
the corresponding constructed
circuit $\Omega$ belongs to $\mathcal{A}^{1,\alpha}(r_1,L_1,\tilde{R})$, for some suitable
positive constants $r_1\leq r$, $L_1\geq L$, $\tilde{R}\geq R$ and $\alpha$, $0<\alpha<1$.
Then we can find positive constants $c_4$, $0<c_4\leq 1$, $C_5$ and $C_6$, depending on
$r_1$, $L_1$, $\tilde{R}$ and $\alpha$ only,
such that if $\tilde{d}_1(\Omega_0,\Omega)\leq c_4r_1$, then
we can subdivide $\partial \Omega$ into smooth curves $J_i$, $i=1\ldots,n$,
which are pairwise internally disjoint, such that for any $i=1,\ldots,n$ we have
$$d_H(J_i,I_i)\leq C_5\tilde{d}_1(\Omega_0,\Omega)$$
and
$$\mathrm{length}(I_i)-2C_5\tilde{d}_1(\Omega_0,\Omega)\leq 
\mathrm{length}(J_i)\leq \mathrm{length}(I_i)+C_6(\tilde{d}_1(\Omega_0,\Omega))^{\alpha/(2\alpha+2)}.$$
Therefore,
$$-2nC_5\tilde{d}_1(\Omega_0,\Omega)\leq P(\Omega)-P(\Omega_0)\leq
nC_6(\tilde{d}_1(\Omega_0,\Omega))^{\alpha/(2\alpha+2)}.$$

By these reasonings it might seem that we may choose
to measure the distance between the desired circuit $\Omega_0$ and the reconstructed one $\Omega$ by using any of these distances. However, there are several reasons to prefer the distance $d_3$, which we actually choose. In fact, 
it is easier to compute than $d_1$ and $\tilde{d}_1$, it can be extended in a natural way from characteristic functions to any $BV$ function by using
$d_{st}$, and should provide a better approximation of the desired circuit than $d_2$, which seems to be too weak for this purpose.

\subsection{Convolutions of characteristic functions and Gaussian distributions}

We recall that
$G(x)=(2\pi)^{-1}\rme ^{-|x|^2/2}$, $x\in \mathbb{R}^2$, and let us note that
$\hat{G}(\xi)=\rme ^{-|\xi|^2/2}$, $\xi\in\mathbb{R}^2$. Moreover, $\|G\|_{L^1(\mathbb{R}^2)}=1$.
For any positive constant $s$ we denote by $G_s(x)=s^{-2}G(x/s)$, $x\in \mathbb{R}^2$.
We note that $\|G_s\|_{L^1(\mathbb{R}^2)}=1$ and $\widehat{G_s}(\xi)=\hat{G}(s\xi)$, $\xi\in\mathbb{R}^2$.

Let $D$ be a bounded open set which is Lipschitz with constants $R_0$ and $L$ and let $\chi_{D}$
be its characteristic function.
We investigate how $\chi_{D}$ is perturbed if we convolute it with $G$. We call $v=\chi_{D}\ast G$,
that is
$$v(x)=\int_{\mathbb{R}^2}\chi_{D}(x-y)G(y)\rmd y=
\int_{\mathbb{R}^2}\chi_{D}(y)G(x-y)\rmd y,\quad x\in\mathbb{R}^2.$$
We recall that we shall use the positive constants $\delta_1$, $\delta_2$ and $m_1$, and the sets $M_{\delta_1}$ and
$N_{\delta}$ introduced at the beginning of Subsection~\ref{geomsubsec}.

\begin{prop}\label{firstapproxprop}
Under the previous notation and assumptions, let
us fix $\delta$, $0<\delta\leq \delta_2/4$. Then
there exist constants $R_0\geq 1$, $\tilde{h}$, $0<\tilde{h}\leq 1/24$, and $a_1>0$,
depending on $L$ and $\delta$ only, such that the following estimates hold.
For any $x\in\mathbb{R}^2$,
\begin{equation}\label{a1}
\text{if }\tilde{h}<v(x)<1-\tilde{h},\text{ then }x\in \overline{B}_{m_1\delta R_0}(\partial D),
\end{equation}
and for any $x\in\partial D$,
\begin{equation}\label{a2}
\text{if }y\in N_{\delta}(x),\text{ then }\nabla v(y)\cdot (-e_2)\geq a_1.
\end{equation}
\end{prop}

\proof{.} If
$x\not \in  \overline{B}_{m_1\delta R_0}(\partial D)$, then we have
$$v(x)\leq \rme^{-m_1^2\delta^2R_0^2/2},\text{ if }x\not\in D,$$
and
$$v(x)\geq 1-\rme^{-m_1^2\delta^2R_0^2/2},\text{ if }x\in D.$$
Consequently, provided $\tilde{h}=\rme^{-m_1^2\delta^2R_0^2/2}\leq 1/24$, we may conclude that
\eqref{a1} holds.

Let us take $x\in \partial D$ and
$y\in N_{\delta}(x)$.
Then, denoting by $\nu$ the exterior unit normal vector to $ D$,
$$\nabla v(y)\cdot (-e_2)=
\int_{\partial D}G(y-z)\nu(z)\cdot e_2\rmd \mathcal{H}^1(z).$$
Therefore,
\begin{multline*}
\nabla v(y)\cdot (-e_2)=
\int_{\partial D\cap \overline{B}_{2\delta R_0}(y)}G(y-z)\nu(z)\cdot e_2\rmd \mathcal{H}^1(z)+\\
\int_{\partial D\cap(B_{R_0/2}(y)\backslash \overline{B}_{2\delta R_0}(y)) }G(y-z)\nu(z)\cdot e_2\rmd \mathcal{H}^1(z)+\\
\int_{\partial D\backslash B_{R_0/2}(y)}G(y-z)\nu(z)\cdot e_2\rmd \mathcal{H}^1(z)=A+B+C.
\end{multline*}

Since $B_{R_0/2}(y)$ is contained in $B_{R_0}(x)$,
for any $z\in \partial D\cap B_{R_0/2}(y)$,
we have $\nu(z)\cdot e_2\geq c_1>0$ where $c_1$ is a constant depending on $L$ only.
Moreover, the length of
$\partial D\cap \overline{B}_{2\delta R_0}(y)$ is also bounded from below by $c_2\delta R_0$, $c_2>0$ depending on $L$ only.
Therefore, we obtain that $A\geq c_1c_2\delta R_0\rme^{-2\delta^2R_0^2}$ and $B\geq 0$.

For what concerns the term $C$, with the help of \eqref{length1}, we can find a constant $C_1$, depending on $L$ only,
such that, for any $R_0\geq 1$, we have
$$|C|\leq C_1R_0\rme^{-R_0^2/8}.$$

Therefore, we can find $R_0\geq 1$, depending on $L$ and $\delta$ only, such that 
$\tilde{h}=\rme^{-m_1^2\delta^2R_0^2/2}\leq1/24$, 
and $2C_1\rme^{-R_0^2/8}\leq c_1c_2\delta \rme^{-2\delta^2R_0^2}$.
We set
$a_1=(1/2)c_1c_2\delta R_0 \rme^{-2\delta^2R_0^2}$ and the proof is concluded.\cvd

\begin{oss}\label{R0oss}
Without loss of generality, we may choose $R_0$ such that it also satisfies
\begin{equation}\label{R0choice}
\|\nabla G\|_{L^1(\mathbb{R}^2\backslash B_{R_0/2})}\leq (1/12)a_1.
\end{equation}
\end{oss}

In the sequel, we shall fix $\delta=\delta_2/4$ and $R_0$ as the corresponding constant in
Proposition~\ref{firstapproxprop} such that \eqref{R0choice} holds. We note that, in this case,
$\delta$ and $R_0$ depend on $L$ only. We shall also fix a constant $R\geq 10R_0$.
We recall that, with a slight abuse of notation, we identify $L^1(B_R)$ with the set of real valued $L^1(\mathbb{R}^2)$ functions that are
equal to zero almost everywhere outside $B_R$.
The same proof of Proposition~\ref{firstapproxprop} allows us to prove
this corollary.

\begin{cor}\label{cors}
For any $s$, $0<s\leq 1$, let $r=sR_0$ and let $D$ be a bounded open set which is Lipschitz with constants
$r$ and $L$. Let $v_s=\chi_{D}\ast G_s$. Then, for any $x\in\mathbb{R}^2$,
$$\text{if }\tilde{h}<v_s(x)<1-\tilde{h},\text{ then }x\in \overline{B}_{m_1\delta r}(\partial  D),$$
and for any $x\in\partial D$,
$$\text{if }y\in N_{\delta}(x),\text{ then }\nabla v_s(y)\cdot (-e_2)\geq a_1R_0/r=a_1/s.$$
\end{cor}

We conclude this part with the following perturbation argument. Let us consider a function $\psi$ such that either
$\psi
\in C^1(\mathbb{R}^2)\cap W^{1,1}(\mathbb{R}^2)$ or $\psi
\in W^{2,1}(\mathbb{R}^2)$ and that, for some $\tilde{\delta}>0$,
$$\|\psi\|_{W^{1,1}(\mathbb{R}^2)}\leq \tilde{\delta}.$$
Let $\tilde{G}=G+\psi$. Then the following result holds.

\begin{cor}\label{cors2}
Let us assume that $\tilde{\delta}\leq \min\{\tilde{h},a_1/2\}$.

For any $s$, $0<s\leq 1$, let $r=sR_0$ and let $D$ be a bounded open set which is Lipschitz with constants
$r$ and $L$. Let $v_s=\chi_{D}\ast \tilde{G}_s$. Then, for any $x\in\mathbb{R}^2$,
\begin{equation}\label{in}
\text{if }2\tilde{h}<v_s(x)<1-2\tilde{h},\text{ then }x\in \overline{B}_{m_1\delta r}(\partial  D),
\end{equation}
and for any $x\in\partial D$,
\begin{equation}\label{gradlower}
\text{if }y\in N_{\delta}(x),\text{ then }\nabla v_s(y)\cdot (-e_2)\geq a_1R_0/2r=a_1/(2s).
\end{equation}
\end{cor}

\proof{.} It follows immediately from the previous corollary and Proposition~\ref{convprop}. We first notice that in either cases $\chi_D\ast \tilde{G}_s\in C^1(\mathbb{R}^2)$. Moreover
we have
$$\|\chi_D\ast \tilde{G}_s-\chi_D\ast G_s\|_{L^{\infty}(\mathbb{R}^2)}\leq \|\chi_D\|_{L^{\infty}(\mathbb{R}^2)}\|\psi_s\|_{L^1}\leq \|\psi\|_{L^1(\mathbb{R}^2)}$$
and
\begin{multline*}
\|\nabla (\chi_D\ast \tilde{G}_s-\chi_D\ast G_s)\|_{L^{\infty}(\mathbb{R}^2)}=
\|\chi_D\ast (\nabla\psi_s)\|_{L^{\infty}(\mathbb{R}^2)}\leq\\ \|\chi_D\|_{L^{\infty}(\mathbb{R}^2)}\|\nabla \psi_s\|_{L^1}\leq \|\nabla \psi\|_{L^1(\mathbb{R}^2)}/s.
\end{multline*}
Thus the conclusion follows.\cvd

\section{Relationship between a mask and its image intensity}\label{dirsec}

In this section we study the relationship between a function representing a mask (not
necessarily a characteristic function of a domain) and its associated image intensity.
We recall the notation used.
We fix $\delta=\delta_2/4$ and $R_0$ as the corresponding constant in
Proposition~\ref{firstapproxprop} such that \eqref{R0choice} holds. We note that, in this case,
$\delta$ and $R_0$ depend on $L$ only. We shall also fix a constant $R\geq 10R_0$.
We recall that, with a slight abuse of notation, we identify $L^1(B_R)$ with the set of real valued $L^1(\mathbb{R}^2)$ functions that are
equal to zero almost everywhere outside $B_R$.
We recall that
$A=\{u\in L^1(B_R)\,:\ 0\leq u\leq 1\text{ a.e. in }B_R\}$.

Fixed $\tilde{\delta}>0$,
we assume that
$\psi\in W^{2,1}(\mathbb{R}^2)$ and that
$$\|\psi\|_{W^{1,1}(\mathbb{R}^2)}\leq \tilde{\delta}.$$
We denote $\tilde{G}=G+\psi$ and, for any $s$, $0<s\leq 1$, we define the operator
$\mathcal{T}_s:L^1(B_R)\to W^{2,1}(\mathbb{R}^2)$ as follows
$$\mathcal{T}_s(u)=u\ast \tilde{G}_s,\quad\text{for any }u\in L^1(B_R).$$
The point spread function we use, $T$, can
be described in general by the function $\tilde{G}$.  Therefore a
study of properties of convolutions with $\tilde{G}$ will be useful.

We remark that the following continuity properties of the operator
$\mathcal{T}_s$ hold. For any $p$, $1\leq p\leq +\infty$, and any $u\in L^1(B_R)$, we have, for an absolute constant $C$, 
\begin{align}\label{c_1}
&\|\mathcal{T}_s(u)\|_{L^p(\mathbb{R}^2)}\leq \|\tilde{G}\|_{L^1(\mathbb{R}^2)}\|u\|_{L^p(\mathbb{R}^2)},\\\label{c_2}
&\|\nabla \mathcal{T}_s(u)\|_{L^p(\mathbb{R}^2)}\leq (C/s)\|\nabla \tilde{G}\|_{L^1(\mathbb{R}^2)}\|u\|_{L^p(\mathbb{R}^2)},\\\label{c_3}
&\|D^2 \mathcal{T}_s(u)\|_{L^p(\mathbb{R}^2)}\leq (C/s^2)\|D^2 \tilde{G}\|_{L^1(\mathbb{R}^2)}\|u\|_{L^p(\mathbb{R}^2)}.
\end{align}

Let $J\in C^{0}(\mathbb{R}^2)$.
For any $u\in L^1(B_R)$, we define $U\in L^1(\mathbb{R}^4)$ as follows
$$U(x,y)=u(x)u(y)J(x-y),\quad\text{for any }x,y\in\mathbb{R}^2.$$
Then, for any  $s$, $0<s\leq 1$, we define 
$H_s\in W^{2,1}(\mathbb{R}^4)$ in the following way
$$H_s(x,y)=\tilde{G}_s(x)\tilde{G}_s(y),\quad\text{for any }x,y\in\mathbb{R}^2.$$
Therefore, for any $p$, $1\leq p\leq +\infty$, and any $u\in L^1(B_R)$, we have, for an absolute constant $C$,
\begin{align}\label{cc_1}
&\|U\ast H_s\|_{L^p(\mathbb{R}^4)}\leq \|\tilde{G}\|^2_{L^1(\mathbb{R}^2)}
\|J\|_{L^{\infty}(B_{2R})}
\|u\|^2_{L^p(\mathbb{R}^2)},\\\label{cc_2}
&\|\nabla(U\ast H_s)\|_{L^p(\mathbb{R}^4)}\leq (C/s)
\|\tilde{G}\|_{L^1(\mathbb{R}^2)}
\|\nabla \tilde{G}\|_{L^1(\mathbb{R}^2)}\|J\|_{L^{\infty}(B_{2R})}\|u\|^2_{L^p(\mathbb{R}^2)},\\\label{cc_3}
&\|D^2 (U\ast H_s)\|_{L^p(\mathbb{R}^4)}\leq (C/s^2)\|\tilde{G}\|^2_{W^{2,1}(\mathbb{R}^2)}\|J\|_{L^{\infty}(B_{2R})}\|u\|^2_{L^p(\mathbb{R}^2)}.
\end{align}

Let us fix $p>4$ and let $\alpha=1-4/p$, $0<\alpha<1$. Then, if $u\in A$ we have
$U\ast H_s\in C^{1,\alpha}(\mathbb{R}^4)$ and, for some absolute constant $C$ depending on $p$,
$$\|U\ast H_s\|_{C^{1,\alpha}(\mathbb{R}^4)}\leq C\|U\ast H_s\|_{W^{2,p}(\mathbb{R}^4)}.$$

We define
$\mathcal{P}_{J,s}:A\to C^{1,\alpha}(\mathbb{R}^2)$ and $\mathcal{P}_{1,s}:A\to C^{1,\alpha}(\mathbb{R}^2)$ as follows.
For any $u\in A$
$$\mathcal{P}_{J,s}(u)(x)=\int_{\mathbb{R}^2} \int_{\mathbb{R}^2} u(\xi)\tilde{G}_s(x-\xi) J(\xi-\eta) \tilde{G}_s(x-\eta) u(\eta) \rmd\xi \rmd\eta,\quad x\in\mathbb{R}^2$$
and
$$\mathcal{P}_{1,s}(u)=(\mathcal{T}_s(u))^2.$$
We notice that the two definitions are consistent when $J\equiv 1$ and that
$$\mathcal{P}_{J,s}(u)(x)=(U\ast H_s)(x,x),\quad x\in\mathbb{R}^2.$$

Putting together the previous estimates we obtain
the following result. We recall that we have fixed a number $p>4$ and that $\alpha=1-4/p$.

\begin{prop}\label{coherenceproposition}
Under the previous notation and assumptions,
let $\varepsilon=\|J-1\|_{L^{\infty}(B_{2R})}$. Then for any $u\in A$ and any $s$, $0<s\leq 1$, we have, for some absolute constant $C$ depending on $p$,
\begin{multline*}
\|\mathcal{P}_{J,s}(u)\|_{C^{0,\alpha}(\mathbb{R}^2)}\leq ((1+\varepsilon)C/s)\|\tilde{G}\|^2_{W^{1,1}(\mathbb{R}^2)}\|u\|^2_{L^{p}(B_R)}\leq\\ ((1+\varepsilon)C/s)\|\tilde{G}\|^2_{W^{1,1}(\mathbb{R}^2)}\|u\|_{L^1(B_R)}^{2/p}.
\end{multline*}
The same estimate holds also for the gradient, namely
\begin{multline*}
\|\nabla \mathcal{P}_{J,s}(u)\|_{C^{0,\alpha}(\mathbb{R}^2)}\leq
((1+\varepsilon)C/s^2)\|\tilde{G}\|^2_{W^{2,1}(\mathbb{R}^2)}\|u\|^2_{L^p(B_R)}\leq\\
((1+\varepsilon)C/s^2)\|\tilde{G}\|^2_{W^{2,1}(\mathbb{R}^2)}\|u\|_{L^1(B_R)}^{2/p}.
\end{multline*}

Furthermore, we have
\begin{equation}\label{0coher}
\|\mathcal{P}_{J,s}(u)-\mathcal{P}_{1,s}(u)\|_{L^{\infty}(\mathbb{R}^2)}=\|\mathcal{P}_{J-1,s}(u)\|_{L^{\infty}(\mathbb{R}^2)}\leq
\|\tilde{G}\|^2_{L^1(\mathbb{R}^2)} \varepsilon
\end{equation}
and, for some absolute constant $C$,
\begin{equation}\label{1coher}
\|\nabla (\mathcal{P}_{J,s}(u)-\mathcal{P}_{1,s}(u))\|_{L^{\infty}(\mathbb{R}^2)}\leq C\|\tilde{G}\|_{L^1(\mathbb{R}^2)}
\|\nabla \tilde{G}\|_{L^1(\mathbb{R}^2)}\varepsilon/s.
\end{equation}
\end{prop}

Although $\mathcal{P}_{J,s}$ is nonlinear in its argument $u$, by a simple adaptation of the previous reasonings, we obtain that for any $u_1$, $u_2\in A$, and for some absolute constant $C$ depending on $p$, we have the following corresponding estimates
\begin{multline*}
\|\mathcal{P}_{J,s}(u_1)-\mathcal{P}_{J,s}(u_2)\|_{C^{0,\alpha}(\mathbb{R}^2)}\leq ((1+\varepsilon)C/s)\|\tilde{G}\|^2_{W^{1,1}(\mathbb{R}^2)}R^{2/p}
\|u_1-u_2\|_{L^{p}(B_R)}
\leq\\ 2((1+\varepsilon)C/s)\|\tilde{G}\|^2_{W^{1,1}(\mathbb{R}^2)}R^{2/p}
\|u_1-u_2\|^{1/p}_{L^1(B_R)}
\end{multline*}
and
\begin{multline*}
\|\nabla (\mathcal{P}_{J,s}(u_1)-\mathcal{P}_{J,s}(u_2))\|_{C^{0,\alpha}(\mathbb{R}^2)}\leq
((1+\varepsilon)C/s^2)\|\tilde{G}\|^2_{W^{2,1}(\mathbb{R}^2)}R^{2/p}\|u_1-u_2\|_{L^{p}(B_R)}\leq\\
2((1+\varepsilon)C/s^2)\|\tilde{G}\|^2_{W^{2,1}(\mathbb{R}^2)}R^{2/p}\|u_1-u_2\|_{L^1(B_R)}^{1/p}.
\end{multline*}
Therefore, $\mathcal{P}_{J,s}:A\to C^{1,\alpha}(\mathbb{R}^2)$
is Lipschitz continuous with respect to the $L^p$ norm and H\"older continuous
with exponent $1/p$ with respect to the $L^1$ norm.

We fix $\tilde{\delta}$ such that $0<\tilde{\delta}\leq\min\{\tilde{h},a_1/2\}$, with $\tilde{h}$, $0<\tilde{h}\leq 1/24$, and $a_1>0$ as in Proposition~\ref{firstapproxprop}, thus depending on $L$ only. We define the corresponding $T$ and $s_0$ as in Lemma~\ref{approxlemma}.
We finally fix $s=(k\text{NA})^{-1}$, $0<s\leq 1/s_0$, and $\sigma$, $0<\sigma\leq s$.
Then we define
$$I_{s,\sigma}(u)=\mathcal{P}_{J,ss_0}(u),\quad\text{for any }u\in A,$$
where $J$ is given by \eqref{Jdef}. We recall that in \eqref{Kdef} we defined
$K=T_{ss_0}$, therefore for any open set $D\subset B_R$ we have that
$I_{s,\sigma}(\chi_D)$ is the light intensity on the image plane corresponding to the mask $D$, see \eqref{Hopkins}.

We denote by $\mathcal{H}:\mathbb{R}\to\mathbb{R}$ the Heaviside function
such that $\mathcal{H}(t)=0$ for any $t\leq 0$ and $\mathcal{H}(t)=1$ for any $t> 0$.
For any constant $h$ we set $\mathcal{H}_h(t)=\mathcal{H}(t-h)$
for any $t\in\mathbb{R}$.
Then, for any $h$, $0<h<1$, any $s$, $0<s\leq 1/s_0$, and any $\sigma$, $0<\sigma\leq s$,
we define the operator
$\mathcal{W}:A\to L^{\infty}(\mathbb{R}^2)$ as follows
\begin{equation}\label{Wdefinition}
\mathcal{W}(u)=\mathcal{H}_h(I_{s,\sigma}(u)),\quad \text{for any }u\in A.
\end{equation}
Clearly, for any $u\in A$, $\mathcal{W}(u)$ is the characteristic function of an open set,
which we shall call $\Omega(u)$. That is
\begin{equation}\label{Omegedefinition}
\Omega(u)=\{x\in\mathbb{R}^2\,:\ I_{s,\sigma}(u)(x)>h\},\quad \text{for any }u\in A.
\end{equation}
In other
words, $\chi_{\Omega(u)}=\mathcal{W}(u)=\mathcal{H}_h(I_{s,\sigma}(u))$.
Moreover, whenever $u=\chi_D$, where $D$ is an open set contained in $B_R$,
we shall denote $\Omega(D)=\Omega(\chi_D)$.

The final, and crucial, result of this section is the following. 

\begin{teo}\label{mainteo}
Let us fix a positive constant $L$. Let $\delta=\delta_2/4$ and let $R_0$ be as in Proposition~\textnormal{\ref{firstapproxprop}} and such that \eqref{R0choice} holds.
Let us also fix $R\geq 10 R_0$ and $p$, $p>4$, and $\alpha=1-4/p$.

We fix $\tilde{\delta}$ such that $0<\tilde{\delta}\leq\min\{\tilde{h},a_1/2\}$, with $\tilde{h}$,  $0<\tilde{h}\leq 1/24$, and $a_1>0$ as in Proposition~\textnormal{\ref{firstapproxprop}}, thus depending on $L$ only. We define the corresponding $T$ and $s_0$ as in Lemma~\textnormal{\ref{approxlemma}}.
We finally fix $s=(k\text{NA})^{-1}$, $0<s\leq 1/s_0$, and $\sigma$, $0<\sigma\leq s$.
Then, for any $u\in A$ we define
$$I_{s,\sigma}(u)=\mathcal{P}_{J,ss_0}(u)$$
where $J$ is given by \eqref{Jdef}.

Then
for any $h$, $1/3\leq h\leq 2/3$, and any $s$, $0<s\leq 1/s_0$,
we can find positive constants  $\tilde{\sigma}_0$, $0<\tilde{\sigma}_0\leq 1$, and $\gamma_0$,
depending on $L$, $R$, $\|D^2(T-G)\|_{L^1(\mathbb{R}^2)}$, $p$,
and $ss_0$ only,
such that for any
$\sigma$, $0<\sigma\leq \tilde{\sigma}_0s$, and any $\gamma$, $0<\gamma\leq \gamma_0$,
the following holds.

Let $\mathcal{A}=\mathcal{A}^{0,1}(r,L,R)$, where $r=ss_0R_0$. Let $\tilde{R}=R+2m_1\delta R_0$, where $m_1$, $0<m_1\leq 1$, depends on $L$ only. 
Then, for any $u\in \mathcal{A}_{\gamma}$, we have that $\Omega(u)\Subset B_{\tilde{R}}$ and
$\Omega(u)\in \mathcal{A}^{1,\alpha}(r_1,L_1,\tilde{R})$.
Here $r_1=ss_0\tilde{R}_0\leq r$, where $\tilde{R}_0 \leq \delta_1R_0/8$ depends on $L$ only,
whereas $L_1\geq L$ depends on
$L$, $R$, $\|D^2(T-G)\|_{L^1(\mathbb{R}^2)}$, $p$ and $ss_0$ only.
 
Moreover, the map $\mathcal{W}:\mathcal{A}_{\gamma}\to BV(B_{\tilde{R}})$ is uniformly continuous with respect to the $L^1$ norm on $\mathcal{A}_{\gamma}$
and the distance $d_{st}$ on $BV(B_{\tilde{R}})$.
\end{teo}

\begin{oss}\label{ossimp}
We observe that the distance $d_{st}$ in $BV(B_{\tilde{R}})$ between $\mathcal{W}(u_1)$
and $\mathcal{W}(u_2)$
corresponds to the distance $d_3$ related to $B_{\tilde{R}}$ between $\Omega(u_1)$ and $\Omega(u_2)$.
\end{oss}

\proof{ of Theorem~\ref{mainteo}.}
The proof is a consequence of the previous analysis.
We fix $s$, $0<s\leq 1/s_0$,
and $h$, $1/3\leq h\leq 2/3$. 

Let us begin with the following preliminary case. Let $u=\chi_D$, where $D\in\mathcal{A}$,
and let $v=\mathcal{T}_{ss_0}(u)$ and
 $\tilde{W}=\mathcal{H}_h((\mathcal{T}_{ss_0}(u))^2)$.
 We apply Corollary~\ref{cors2} and we obtain the following results.

If $\tilde{\Omega}$ is the open set such that $\tilde{W}=\chi_{\tilde{\Omega}}$, then, by \eqref{in}, we notice that $\partial
\tilde{\Omega}\subset \overline{B}_{m_1\delta r}(\partial D)$ and that $(D\backslash \overline{B}_{m_1\delta r}(\partial D))\subset \tilde{\Omega}$ and
$(\mathbb{R}^2\backslash \overline{B}_{m_1\delta r}(\overline{D}))\cap \tilde{\Omega}=\emptyset$. Therefore
 $\tilde{\Omega}\Subset B_{\tilde{R}}$.

We take any
$x\in \partial D$ and any $y\in M_{\delta_1/2}(x)$, with respect to the coordinate system depending on $x$. Then we consider
the points $y^-=y-\delta re_2$ and $y^+=y+\delta re_2$. We have that
$y^{\pm}\in\partial N_{\delta}(x)\backslash \overline{B}_{m_1\delta r}(\partial D)$.  Moreover, 
$y^-\in D$ and $v(y^-)\geq 11/12$, whereas $y^+\not\in D$ and
$v(y^+)\leq 1/12$. Let us call $\tilde{y}^+=y+t_0\delta re_2$, where $t_0\in (-1,1]$, $t_0$ depends on $y$,
and $v(\tilde{y}^+)=1/12$ whereas $v(y+t\delta re_2)< 1/12$ for any $t\in (t_0,1]$.
Then we use \eqref{gradlower} and we obtain that, for any $t\in [-1,t_0]$,
$v(y+t\delta re_2)\geq 1/12$ and
$$-\frac{\rmd}{\rmd t}(v(y+t\delta re_2))^2\geq \delta ra_1/(12ss_0).$$
We may conclude that there exists a function $\varphi_{1}:[-\delta_1r/2,\delta_1r/2]\to \mathbb{R}$
such that, for any $y=(y_1,y_2)\in N_{\delta}(x)$ with $|y_1|\leq \delta_1r/2$,
$(v(y))^2=h$ if and only if
$y_2=\varphi_1(y_1)$.
We recall that
\begin{equation}\label{aa1}
\|v\|_{L^{\infty}(\mathbb{R}^2)}\leq C_1,\quad\|\nabla v\|_{L^{\infty}(\mathbb{R}^2)}\leq C_1/(ss_0),
\end{equation}
where $C_1$ is an absolute constant, and
\begin{equation}\label{aa2}
\|\nabla v\|_{C^{0,\alpha}(\mathbb{R}^2)}\leq C_2/(ss_0)^2,
\end{equation}
where $C_2$ depends on $R$, $\|D^2(T-G)\|_{L^1(\mathbb{R}^2)}$ and $p$ only.

We obtain that $v^2$ is a $C^{1,\alpha}$ function and, by the implicit function theorem, we conclude
that the function $\varphi_1$ is actually $C^{1,\alpha}$.
We observe that
$$\|\varphi'_1\|_{L^{\infty}[-\delta_1r/2,\delta_1r/2]}\leq C_3,$$
where $C_3$ depends on $L$ only. Without loss of generality, by a translation
we may assume that
$\varphi_1(0)=0$, thus
$\|\varphi_1\|_{L^{\infty}[-\delta_1r/2,\delta_1r/2]}\leq C_3\delta_1r/2$.
Finally, for any $t_1$, $t_2\in [-\delta_1r/2,\delta_1r/2]$,
$$|\varphi'_1(t_1)-\varphi'_1(t_2)|\leq (C_4/(ss_0))|t_1-t_2|^{\alpha},$$
where $C_4$ is a constant depending on $L$, $R$, $\|D^2(T-G)\|_{L^1(\mathbb{R}^2)}$ and $p$ only.

Then, it is not difficult to prove that for some $r_1=ss_0\tilde{R}_0\leq r$, with $\tilde{R}_0 \leq \delta_1R_0/8$ depending on $L$ only, we can find
$L_1\geq L$, depending on $L$, $R$, $\|D^2(T-G)\|_{L^1(\mathbb{R}^2)}$, $p$ and $ss_0$ only,
such that $\tilde{\Omega}\in \mathcal{A}^{1,\alpha}(r_1,L_1,\tilde{R})$.
Let us also remark that we have obtained that
$\tilde{d}_1(\tilde{\Omega},D)\leq\delta r$.

Let us call $\varepsilon=\varepsilon(s,\sigma)=
\varepsilon(\sigma/s)=
\|J-1\|_{L^{\infty}(B_{2R})}$.
We notice that, as $\sigma/s\to 0^+$, we have that $\varepsilon$ goes to $0$ as well.
We also assume, without loss of generality, that $\varepsilon$ is increasing with respect to the variable $\sigma/s$.
Let us recall that, for any $u\in A$, if $w=I_{s,\sigma}(u)$, with
$0<s\leq 1/s_0$ and $0<\sigma\leq s$,
then
\begin{equation}\label{b1}
\|w\|_{L^{\infty}(\mathbb{R}^2)}\leq C_5,\quad \|\nabla w\|_{L^{\infty}(\mathbb{R}^2)}\leq C_5/(ss_0),
\end{equation}
and
\begin{equation}\label{b2}
\|\nabla w\|_{C^{0,\alpha}(\mathbb{R}^2)}\leq C_6/(ss_0)^2,
\end{equation}
where $C_5$ is an absolute constant
and $C_6$ depends on $R$, $\|D^2(T-G)\|_{L^1(\mathbb{R}^2)}$ and $p$ only.

For positive constants $\tilde{\sigma}_0$, $0<\tilde{\sigma}_0\leq 1$, and $\gamma_0$, to be precised later,
let us fix $\sigma$, $0<\sigma\leq \tilde{\sigma}_0s$ and $\gamma$, $0<\gamma\leq\gamma_0$.
We take $u\in\mathcal{A}_{\gamma}$, $v=\mathcal{T}_{ss_0}(u)$,
$w=I_{s,\sigma}(u)$, and $D\in\mathcal{A}$ such that $\|u-\chi_D\|_{L^1(B_R)}\leq \gamma$. Then we use Proposition~\ref{coherenceproposition} to infer that
\begin{multline*}
\left\|I_{\sigma,s}(u)-(\mathcal{T}_{ss_0}(\chi_D))^2\right\|_{L^{\infty}(\mathbb{R}^2)}\leq
\|I_{s,\sigma}(u)-v^2\|_{L^{\infty}(\mathbb{R}^2)}+
\|(\mathcal{T}_{ss_0}(u))^2-(\mathcal{T}_{ss_0}(\chi_D))^2\|_{L^{\infty}(\mathbb{R}^2)}
\leq \\
\|\tilde{G}\|^2_{L^1(\mathbb{R}^2)}\varepsilon+(2C/(ss_0))\|\tilde{G}\|_{L^1(\mathbb{R}^2)}\|\tilde{G}\|_{W^{1,1}(\mathbb{R}^2)}\|u-\chi_D\|_{L^p(B_R)}\leq
(C_7/(ss_0))\left(\varepsilon+\gamma^{1/p}\right),
\end{multline*}
where $C$ is an absolute constant
and consequently $C_7$ depends on $p$ only.

Analogously, we can prove that
$$\left\|\nabla\left(I_{s,\sigma}(u)-(\mathcal{T}_{ss_0}(\chi_D))^2\right)\right\|_{L^{\infty}(\mathbb{R}^2)}\leq (C_8/(ss_0)^2)\left(\varepsilon+\gamma^{1/p}\right),$$
where the constant $C_8$ depends on $\|D^2(T-G)\|_{L^1(\mathbb{R}^2)}$ and $p$ only.

We now choose the positive constants $\tilde{\sigma}_0$ and $\gamma_0$ in such a way that
$$2(C_7/(ss_0))\left(\varepsilon(\tilde{\sigma}_0)+\gamma_0^{1/p}\right)\leq 1/6$$
and
$$(C_8/(ss_0))\left(\varepsilon(\tilde{\sigma}_0)+\gamma_0^{1/p}\right)\leq a_1/24.$$
Clearly, $\tilde{\sigma}_0$ and $\gamma_0$ depends on $L$, $\|D^2(T-G)\|_{L^1(\mathbb{R}^2)}$, $p$ and
$ss_0$ only.

Then we can apply to $w=I_{s,\sigma}(u)$ and $\Omega=\Omega(u)$
the same analysis we have used for $v$ and $\tilde{\Omega}$ in the first part of this proof. We may therefore conclude that if $u\in\mathcal{A}_{\gamma}$ and $D\in\mathcal{A}$ is such that $\|u-\chi_D\|_{L^1(B_R)}\leq\gamma$, then
$\Omega\Subset B_{\tilde{R}}$,
$\tilde{d}_1(\Omega,D)\leq \delta r$
and,
taken $r_1$ as before, possibly with a smaller $\tilde{R}_0$ still depending on $L$ only,  we can find
$L_1\geq L$, depending on $L$, $R$, $\|D^2(T-G)\|_{L^1(\mathbb{R}^2)}$, $p$ and $ss_0$ only,
such that $\Omega\in \mathcal{A}^{1,\alpha}(r_1,L_1,\tilde{R})$.

This kind of argument leads us also to show that $\Omega$ shares the same topological properties of $D$, that is for example $\Omega$ and $\partial \Omega$
have the same number of connected components of $D$ and $\partial D$, respectively.

It remains to show the uniform continuity property. We recall that the operator $\mathcal{P}_{J,ss_0}$ is H\"older continuous from $A$, with the $L^1(B_R)$ norm,
into $C^{1,\alpha}(\mathbb{R}^2)$, with its usual norm. This means that there exists a constant $\tilde{C}$ such that for any
$u_1$ and $u_2\in A$, if we call $w_1=I_{s,\sigma}(u_1)$
and $w_2=I_{s,\sigma}(u_2)$, then
$$\|w_1-w_2\|_{C^{1,\alpha}(\mathbb{R}^2)}\leq \tilde{C}\|u_1-u_2\|_{L^1(B_R)}^{1/p}.$$
A simple application of the previous analysis allows us to prove this claim

\begin{claim}\label{claim1}
There exists a function $g:[0,+\infty)\to[0,+\infty)$, which is continuous, increasing and such that $g(0)=0$,
satisfying the following property.
For any $u\in\mathcal{A}_{\gamma}$,  
for any $\varepsilon>0$ and any $x\in\mathbb{R}^2$ we have
\begin{equation}\label{claim}
\text{if }x\not\in B_{\varepsilon}(\partial\Omega(u)),\text{ then }|I_{s,\sigma}(u)-h|> g(\varepsilon).
\end{equation}
\end{claim}

\bigskip

Let us now assume that 
$u_1$ and $u_2$ belong to $\mathcal{A}_{\gamma}$ and let us fix $\varepsilon>0$. We can find $\eta>0$
such that if $\|u_1-u_2\|_{L^1(B_R)}\leq \eta$, then $\|w_1-w_2\|_{L^{\infty}}(\mathbb{R}^2)\leq g(\varepsilon)$.

Let us now take $x\in\partial \Omega(u_1)$, that is $x\in\mathbb{R}^2$ such that
$w_1(x)=h$. We infer that $|w_2(x)-h|\leq g(\varepsilon)$, therefore by the claim
we deduce that $x\in B_{\varepsilon}(\partial\Omega(u_2))$. That is
$\partial \Omega(u_1)\subset B_{\varepsilon}(\partial\Omega(u_2))$. By symmetry,
we conclude that
$\tilde{d}_1(\Omega(u_1),\Omega(u_2))\leq \varepsilon$.
In other words, the map which to any $u\in\mathcal{A}_{\gamma}$ associates the open
set $\Omega(u)$ is uniformly continuous with respect to the $L^1$ norm on $\mathcal{A}_{\gamma}$ and the distance $\tilde{d}_1$. However, we have shown in Subsection~\ref{geomsubsec} that the distances $d_1$, $\tilde{d}_1$, $d_2$ and $d_3$
are topologically equivalent on $\mathcal{A}^{1,\alpha}(r_1,L_1,\tilde{R})$, to which all
$\Omega(u)$ belongs, for any $u\in\mathcal{A}_{\gamma}$. Therefore the
map $\mathcal{A}_{\gamma}\ni u\to \Omega(u)$ is uniformly continuous with respect
to the $L^1$ norm on $\mathcal{A}_{\gamma}$ and any of the distances
$d_1$, $\tilde{d}_1$, $d_2$ and $d_3$ related to $B_{\tilde{R}}$.

We observe that
$$d_2(\Omega(u_1),\Omega(u_2))=
\|\mathcal{W}(u_1)-\mathcal{W}(u_2)\|_{L^1(B_{\tilde{R}+1})}=\|\mathcal{W}(u_1)-\mathcal{W}(u_2)\|_{L^1(B_{\tilde{R}})}$$
whereas 
$$d_3(\Omega(u_1),\Omega(u_2))=
d_2(\Omega(u_1),\Omega(u_2))+|P(\Omega(u_1))-P(\Omega(u_2))|
=d_{st}(\mathcal{W}(u_1),\mathcal{W}(u_2))$$
where $d_{st}$ is here the distance inducing strict convergence in $BV(B_{\tilde{R}})$. Therefore we conclude that $\mathcal{W}:\mathcal{A}_{\gamma}\to
BV(B_{\tilde{R}})$ is uniformly continuous with respect to the
$L^1$ norm on $\mathcal{A}_{\gamma}$ and, on $BV(B_{\tilde{R}})$,
with respect either to the $L^1$ norm or to the $d_{st}$ distance.\cvd

\begin{oss} Let us finally remark that if, instead of taking $h\in [1/3,2/3]$, we simply assume
$0<h<1$, then the same analysis may still be carried over. Clearly we need to change the values of $R_0$ and $\tilde{h}_1$ in Proposition~\ref{firstapproxprop}, so that they depend on $h$ as well. As a consequence also the quantities introduced in the above Theorem~\ref{mainteo} would depend on $h$.
\end{oss}

\section{Analysis of the inverse problem}\label{approxsec}

Throughout this section, we shall keep the notation of Theorem~\ref{mainteo} and we shall also assume that the hypotheses of Theorem~\ref{mainteo}
are satisfied. We shall fix $h$, $1/3\leq h\leq 2/3$, and $s$, $0<s\leq 1/s_0$,
and we shall take $\sigma$, $0<\sigma\leq \tilde{\sigma}_0s$, and $\gamma$, $0<\gamma\leq\gamma_0$,
$\tilde{\sigma}_0$ and $\gamma_0$ as in Theorem~\ref{mainteo}.

We call $\Omega_0$
the circuit to be reconstructed and we shall assume that it  belongs to
$\mathcal{A}=\mathcal{A}^{0,1}(r,L,R)$.

We recall that, by Proposition~\ref{compactprop},
$\mathcal{A}$ is compact with respect to the $d_2$ distance, which corresponds to the
distance induced by the $L^1$ norm for the corresponding
characteristic functions.
Then it is an immediate consequence of the last part of Theorem~\ref{mainteo}, see also
Remark~\ref{ossimp}, that the problem
$$\min_{D\in\mathcal{A}}d_3(\Omega(D),\Omega_0)$$
admits a solution. We note that $\Omega(D)=\Omega(\chi_D)$ and that here $d_3$ is the distance defined in \eqref{d3def} related to $B_{\tilde{R}}$ .

From a numerical point of view, the class $\mathcal{A}$ is rather difficult to handle.
We try to reduce this difficulty by enlarging the class $\mathcal{A}$
to a class of characteristic functions of sets with finite perimeter. In order to keep the lower semicontinuity of the functional, we restrict ourselves to
characteristic functions of sets with finite perimeter which are contained in
$\mathcal{A}_{\gamma}$. Namely, we define the following functional
$F_0:A\to [0,+\infty]$ such that for any $u\in A$ we have
\begin{equation}\label{F0}
F_0(u)=
d_{st}(\mathcal{W}(u),\chi_{\Omega_0})+bP(u),
\end{equation}
where $P$ is the functional defined in \eqref{Pdef}
with $\mcd$ chosen to be $B_R$, $b$ is a positive parameter and
$d_{st}$ is the strict convergence distance in $BV(B_{\tilde{R}})$. We recall that,
whenever $u\in\{0,1\}$ almost everywhere in $B_R$ and $u\in BV(B_{R+1})$, then
$P(u)=P(u,B_{R+1})=|Du|(B_{R+1})$.
Otherwise, $P(u)$, and consequently also $F_0(u)$, is equal to $+\infty$.
Moreover, if $u\in \mathcal{A}_{\gamma}$, in particular if $u=\chi_D$ for some $D\in\mathcal{A}$,  then
$d_{st}(\mathcal{W}(u),\chi_{\Omega_0})=d_3(\Omega(u),\Omega_0)$,
where again $d_3$ is the distance defined in \eqref{d3def} related to $B_{\tilde{R}}$ .

We look for the solution to the following minimization problem
\begin{equation}\label{min0}
\min\{F_0(u)\,:\ u\in \mathcal{A}_{\gamma}\}.
\end{equation}
We notice that such a minimization problem admits a solution.

Even if the class $\mathcal{A}_{\gamma}$ might still be not very satisfactory to handle from a numerical point of view, since it somehow involves handling the class $\mathcal{A}$,
we believe that from a practical point of view such a restriction might be dropped and we might use the class $A\subset L^1(B_R)$ instead. In fact, we have a good initial guess,
given by the target circuit $\chi_{\Omega_0}$, and it is reasonable to assume that the 
optimal mask will be a rather small perturbation of $\Omega_0$ itself. In fact, under our assumptions, by the arguments developed in the proof of Theorem~\ref{mainteo},
we can show that $\Omega(u)$ has the same topological properties of $D$, where $\chi_D$ is the element of $\mathcal{A}$ which is closest to $u$. Therefore if we look for
a set $\Omega(u)$ as close as possible to $\Omega_0$, then at least we need to require that the set $D$ has the same topological properties of $\Omega_0$.
For this reason and
since $\Omega_0\in \mathcal{A}$, it might be essentially the same to perform the minimization in a small neighbourhood of $\mathcal{A}$ or in the whole $A$. On the other hand, again by our assumptions, we notice that whenever the boundary of $\Omega_0$  presents a corner, and this is often case, as $\partial\Omega_0$ is often 
the union of a finite number of segments, then $\Omega_0$ cannot be reconstructed in
an exact way, since $\Omega(u)$, for any $u\in\mathcal{A}_{\gamma}$, is a $C^{1,\alpha}$ set,
thus its boundary cannot have any corner. 

Besides dealing with the class $\mathcal{A}_{\gamma}$,
there are several other difficulties. In particular, computing $F_0(\chi_E)$ for some $E\subset B_R$ is not an easy task, since it involves at least the computation of the perimeters of $E$ and of $\Omega(\chi_E)$. Furthermore, solving a minimization problem in the class of sets of finite perimeter
is not a 
straightforward task from the numerical point of view.

In order to solve these difficulties, we use the following strategy. We approximate, in the sense of $\Gamma$-convergence, the functional 
$F_0$
with a family of functional $\{F_{\varepsilon}\}_{\varepsilon>0}$ which are easier
to compute numerically and are
defined on a set of smooth functions.

As in Section 2.3, we take a $C^{\infty}$ function $\phi:\mathbb{R}\to\mathbb{R}$ such that $\phi$
is nondecreasing, $\phi(t)=0$ for any $t\leq -1/2$ and $\phi(t)=1$ for any $t\geq 1/2$.
For any $\eta>0$ and any $\tau\in\mathbb{R}$, let
$$\phi_{\eta, \tau}(t)=\phi\left(\frac{t-\tau}{\eta}\right),\quad\text{for any }t\in\mathbb{R}.$$
Then we have the following result.

\begin{prop}\label{uniformcontprop}
For any $\eta>0$, let
$\Phi_{\eta}:A\to C^{1,\alpha}(\mathbb{R}^2)$
be defined as
$$\Phi_{\eta}(u)=\phi_{\eta,h}(I_{s,\sigma}(u)),\quad\text{for any }u\in A.$$

Then, for any $\eta$, $0<\eta\leq h$,
$\Phi_{\eta}$ is H\"older continuous, with exponent $1/p$,
from $A$, with the $L^1(B_R)$ norm,
into $C^{1,\alpha}(\mathbb{R}^2)$, with its usual norm.

Furthermore, as $\eta\to 0^+$, $(\mathcal{W}-\Phi_{\eta}):\mathcal{A}_{\gamma}\to BV(B_{\tilde{R}})$ converges uniformly to zero on $\mathcal{A}_{\gamma}$
with respect to the distance $d_{st}$ on $BV(B_{\tilde{R}})$.
\end{prop}

\proof{.} The continuity property of $\Phi_{\eta}$ immediately follows
by the continuity of $\mathcal{P}_{J,ss_0}$ and by the properties of $\phi$. We just note that
the H\"older exponent is fixed, whereas the H\"older constant might depend upon
$\eta$.

About the convergence result, 
we begin by recalling that $\mathcal{W}(u)=\mathcal{H}_h(I_{s,\sigma}(u))$,
$u\in\mathcal{A}_{\gamma}$.
We use Claim~\ref{claim1} introduced in the proof of Theorem~\ref{mainteo}.
We call $t_0=\min(1,\sup\{g(s)\,:\ s\in[0,+\infty)\})$ and $s_0$ the positive real number
 such that $g(s_0)=t_0/2$. We call $g^{-1}:[0,t_0/2]\to [0,s_0]$ the continuous, increasing function which is the inverse of $g$ on such intervals.
For any $\eta$, $0<\eta\leq t_0$,
we infer that $\Phi_{\eta}(u)(x)$ might be different from
$\mathcal{W}(u)(x)$ only if $x\in B_{g^{-1}(\eta/2)}(\partial\Omega(u))$.
By estimates like \eqref{length} and \eqref{neigh}, which are independent of $u\in\mathcal{A}_{\gamma}$, we obtain that
$\|(\mathcal{W}-\Phi_{\eta})(u)\|_{L^1(B_{\tilde{R}})}$ converges to zero, as 
$\eta\to 0^+$, uniformly for $u\in\mathcal{A}_{\gamma}$. 

For any $t\in\mathbb{R}$ and any $u\in\mathcal{A}_{\gamma}$,
we call
$$P(u,t)=P(\{x\in\mathbb{R}^2\,:\ I_{s,\sigma}(u)(x)>t\},B_{\tilde{R}}).$$
It remains to prove that, as $\eta\to 0^+$, $|D(\Phi_{\eta}(u))|(B_{\tilde{R}})=
\int_{B_{\tilde{R}}}|\nabla(\Phi_{\eta}(u))|$
converges to $|D(\mathcal{W}(u))|(B_{\tilde{R}})=P(u,h)$ uniformly for
$u\in\mathcal{A}_{\gamma}$.
We argue in the following way. We have that, for any $\eta$, $0<\eta\leq h$, 
$$|D(\Phi_{\eta}(u))|(B_{\tilde{R}})=
\int_{B_{\tilde{R}}}|\nabla(\Phi_{\eta}(u))|=\int_{B_{\tilde{R}}}|\phi'_{\eta,h}(I_{s,\sigma}(u))||\nabla (I_{s,\sigma}(u))|.$$
Since $\phi'_{\eta,h}\geq 0$, and for $\eta$ small enough, uniformly with respect to $u\in\mathcal{A}_{\gamma}$,
$\phi'_{\eta,h}(I_{s,\sigma}(u))=0$ outside $B_{\tilde{R}}$,
without loss of generality, we have that
$$|D(\Phi_{\eta}(u))|(B_{\tilde{R}})=\int_{\mathbb{R}^2}\phi'_{\eta,h}(I_{s,\sigma}(u))|\nabla (I_{s,\sigma}(u))|.$$
By the coarea formula, 
$$
|D(\Phi_{\eta}(u))|(B_{\tilde{R}})=
\int_{-\infty}^{+\infty}\left(\int_{\{I_{s,\sigma}(u)=t\}}
\phi'_{\eta,h}(t)
\rmd \mathcal{H}^1(y)\right)\rmd t .
$$
Therefore,
$$
|D(\Phi_{\eta}(u))|(B_{\tilde{R}})=
\frac{1}{\eta}\int_{-\infty}^{+\infty}
\phi'\left(\frac{(t-h)}{\eta}\right)P(u,t)
\rmd t =
\int_{-1/2}^{+1/2}
\phi'(s)P(u,s\eta+h)
\rmd s .
$$

Since $|D(\mathcal{W}(u))|(B_{\tilde{R}})=P(u,h)$ and
$\int_{-1/2}^{+1/2}
\phi'(s)\rmd s=1$, we obtain that
$$
\big||D(\Phi_{\eta}(u))|(B_{\tilde{R}})-|D(\mathcal{W}(u))|(B_{\tilde{R}})\big|\leq
\int_{-1/2}^{+1/2}
\phi'(s)|P(u,s\eta+h)-P(u,h)|
\rmd s.
$$

It remains to show that, as $\eta\to 0^+$,
$\sup\{|P(u,t+h)-P(u,h)|\,:\ t\in [-\eta/2,+\eta/2]\}$ goes to zero uniformly with respect to $u\in\mathcal{A}_{\gamma}$. Therefore the proof is concluded by using the following claim.

\begin{claim}\label{claim2}
There exist a positive constant $\eta_0$ and a continuous, increasing function $g_1:[0,\eta_0]\to[0,+\infty)$, such that $g_1(0)=0$, such that
for any $\eta$, $0<\eta\leq \eta_0$, and any $u\in \mathcal{A}_{\gamma}$, we have that
$$\sup\{|P(u,t+h)-P(u,h)|\,:\ t\in [-\eta/2,+\eta/2]\}\leq g_1(\eta).$$
\end{claim}

\bigskip

The proof of Claim~\ref{claim2} is a straightforward, although maybe lengthy to describe,
consequence of the analysis developed in the proof of Theorem~\ref{mainteo}. We leave the details to the reader. We just notice that Claim~\ref{claim2} is a sort of generalization of Claim~\ref{claim1} and the arguments used to prove the two claims are essentially
analogous.\cvd

\bigskip

We are now in the position of describing the approximating functionals and
proving the $\Gamma$-convergence result. Let us a fix a constant $p_1$, $1<p_1<+\infty$, and a continuous function $W:\mathbb{R}\to[0,+\infty)$ such that
$W(t)=0$ if and only if $t\in\{0,1\}$.
Let us denote by $P_{\varepsilon}$, $\varepsilon>0$, the functional defined in
\eqref{modmordef} with $p=p_1$, the function $W$ and  $\mcd=B_R$. We recall that the functional $P$
is defined in \eqref{Pdef}, again with $\mcd=B_R$.

Then, for any $\varepsilon>0$, let us define
$F_{\varepsilon}:A\to [0,+\infty]$
such that for any $u\in A$ we have
\begin{equation}
F_{\varepsilon}(u)=d_{st}(\Phi_{\eta(\varepsilon)}(u),\chi_{\Omega_0})+
bP_{\varepsilon}(u)
\end{equation}
where $\eta:[0,+\infty)\to[0,+\infty)$ is a continuous, increasing function such that
$\eta(0)=0$.

By the direct method, we can prove that each of the functionals
$F_{\varepsilon}$, $\varepsilon>0$, admits a minimum either over $A$ or over
$\mathcal{A}_{\gamma}$.

The $\Gamma$-convergence result is the following.

\begin{teo}\label{gammaconvteo}
Let us consider the metric space $(X,d)$ where $X=\mathcal{A}_{\gamma}$ and
$d$ is the metric induced by the $L^1$ norm. Then, as $\varepsilon\to 0^+$,
$F_{\varepsilon}$ $\Gamma$-converges to $F_0$
on $X$ with respect to the distance $d$.
\end{teo}

\proof{.} Let us fix a sequence $\{\varepsilon_n\}_{n=1}^{\infty}$ of positive numbers converging to zero as $n\to\infty$. Let, for any $n\in\mathbb{N}$, $F_n=F_{\varepsilon_n}$. We need to prove that $\Gamma\textrm{-}\!\lim_nF_n=F$.
Let us also remark that we may extend $F_n$ and $F$ over $L^1(\mathbb{R}^2)$
by setting them equal to $+\infty$ outside $\mathcal{A}_{\gamma}$.
Let us define $\tilde{P}_{\varepsilon}$, $\varepsilon>0$, and $\tilde{P}$
as the functionals which are equal to the functionals $P_{\varepsilon}$ and $P$, respectively,
on $\mathcal{A}_{\gamma}$ and $+\infty$ elsewhere.
We recall that $P_{\varepsilon}$, $\varepsilon>0$, and $P$ are defined in \eqref{modmordef} and in \eqref{Pdef}, respectively, with $p=p_1$ and $\mcd=B_R$. 

We observe that, as a
consequence of Proposition~\ref{uniformcontprop} and of the stability of
$\Gamma$-convergence under uniformly converging continuous perturbations,
it is enough to show that $\Gamma\textrm{-}\!\lim_n\tilde{P}_n=\tilde{P}$,
where $\tilde{P}_n=\tilde{P}_{\varepsilon_n}$, $n\in\mathbb{N}$.
Let us prove this $\Gamma$-convergence result.

The $\Gamma$-liminf inequality is an
immediate consequence of Theorem~\ref{Mod-Morteo} and of
the fact that $\mathcal{A}_{\gamma}$ is a closed subset of $L^1(B_R)$.

For what concerns the recovery sequence, then we argue in the following way.
If $u\in A$ is such that $\|u-\chi_D\|_{L^1(B_R)}<\gamma$, for some $D\in\mathcal{A}$, then we again use Theorem~\ref{Mod-Morteo} to construct a recovery sequence for such
a function $u$, that is a sequence $\{u_n\}_{n=1}^{\infty}$ contained in $\mathcal{A}_{\gamma}$ such that, as $n\to\infty$, $u_n\to u$ in $L^1(B_R)$
and $\tilde{P}_n(u_n)\to \tilde{P}(u)$.

It remains to study the case when $u\in \partial\mathcal{A}_{\gamma}$ and $\tilde{P}(u)<+\infty$. In this case, we have
that $u=\chi_E$, where $E\subset B_R$ is a set of finite perimeter,
and we
pick $D\in\mathcal{A}$ such that
$\|\chi_E-\chi_D\|_{L^1(B_R)}=|E\Delta D|=\gamma$. Then at least one of these two cases must be satisfied. Either there exists $x\in B_R\backslash\overline{D}$
such that
$$\lim_{\rho\to 0^+}\frac{|E\cap B_{\rho}(x)|}{|B_{\rho}(x)|}=1$$
or there exists $x\in D$ such that
$$\lim_{\rho\to 0^+}\frac{|E\cap B_{\rho}(x)|}{|B_{\rho}(x)|}=0.$$
We choose an arbitrary sequence $\{\rho_j\}_{j=1}^{\infty}$ of positive numbers
such that $\lim_j\rho_j=0$.
In the first case, for any $j\in\mathbb{N}$, we choose $E_j$ such that
$\chi_{E_j}=\chi_E(1-\chi_{B_{\rho_j}(x)})$.
In the second case, we choose $E_j$ such that
$\chi_{E_j}=\chi_E(1-\chi_{B_{\rho_j}(x)})+\chi_{B_{\rho_j}(x)}$.
We notice that, in either cases, for any $j\in\mathbb{N}$, $E_j$ is a set of finite perimeter
such that $\|\chi_{E_j}-\chi_D\|_{L^1(B_R)}<\gamma$. Furthermore,
as $j\to\infty$ we have that $\chi_{E_j}\to\chi_E$ in $L^1(B_R)$ and
$P(E_j)\to P(E)$, that is $\tilde{P}(\chi_{E_j})\to \tilde{P}(\chi_E)$. Then the proof may be concluded
by following the arguments of Section~4.2 in \cite{Bra} which we have briefly recalled in the proof of Theorem~\ref{Mod-Morteo}.\cvd

\bigskip

We remark that $\Omega_0\in\mathcal{A}$, therefore we may find a
family $\{\tilde{u}_{\varepsilon}\}_{\varepsilon>0}$ such that,
as $\varepsilon\to 0^+$, $\tilde{u}_{\varepsilon}\to \chi_{\Omega_0}$ in $L^1(B_R)$
and
$P_{\varepsilon}(\tilde{u}_{\varepsilon})\to P(\Omega_0)$.
Without loss of generality, we may assume that, for any $\varepsilon>0$,
$0\leq \tilde{u}_{\varepsilon}\leq 1$ almost everywhere in $B_R$ and that
$\tilde{u}_{\varepsilon}\in\mathcal{A}_{\gamma}$. By Proposition~\ref{uniformcontprop},
we conclude that $F_{\varepsilon}(\tilde{u}_{\varepsilon})\to F_0(\Omega_0)<+\infty$.
We obtain that for any $\varepsilon_0>0$
there exists a constant $C_1$ such that
\begin{equation}\label{uniformbound}
\min_{\mathcal{A}_{\gamma}}F_{\varepsilon}\leq C_1\quad\text{for any }\varepsilon,\
0<\varepsilon\leq \varepsilon_0.
\end{equation}
Obviously, the same property is shared by the minimum values of $F_{\varepsilon}$
over $A$.

It remains to prove that the functionals $F_{\varepsilon}$ are equicoercive
over $\mathcal{A}_{\gamma}$, that is that the following result holds.

\begin{prop}\label{equicoerciveprop}
For any $\varepsilon_0>0$,
there exists a compact subset $\mathcal{K}$ of $\mathcal{A}_{\gamma}$ such that for any $\varepsilon$, $0<\varepsilon\leq\varepsilon_0$, we have
$$\min_{\mathcal{K}}F_{\varepsilon}=\min_{\mathcal{A}_{\gamma}}F_{\varepsilon}.$$
\end{prop}

\proof{.}
Let us take the constant $C_1$ as in \eqref{uniformbound}. 
Let $u_{\varepsilon}\in\mathcal{A}_{\gamma}$, $0<\varepsilon\leq \varepsilon_0$,
be such that $F_{\varepsilon}(u_{\varepsilon})=\min_{\mathcal{A}_{\gamma}}F_{\varepsilon}$.
Then we observe that the set $\{u_{\varepsilon}\}_{0<\varepsilon\leq\varepsilon_0}$
satisfies the properties of Remark~\ref{compactnessoss} for some constant $C$.
Therefore $\{u_{\varepsilon}\}_{0<\varepsilon\leq\varepsilon_0}$ is precompact in 
$L^1(B_R)$ and the proof is concluded.\cvd

\begin{oss}
With an analogous proof, the same result of Proposition~\ref{equicoerciveprop} holds
if we replace $\mathcal{A}_{\gamma}$ with $A$.
\end{oss}

\bigskip

By Theorem~\ref{gammaconvteo} and Proposition~\ref{equicoerciveprop}, we can apply the Fundamental Theorem of $\Gamma$-convergence to conclude with the following result.

\begin{teo}\label{finalteo}
We have that $F_0$ admits a minimum over $\mathcal{A}_{\gamma}$ and 
$$
\min_{\mathcal{A}_{\gamma}} F_0=
\lim_{\varepsilon\to 0^+}\inf_{\mathcal{A}_{\gamma}}
F_{\varepsilon}=\lim_{\varepsilon\to 0^+}\min_{\mathcal{A}_{\gamma}}
F_{\varepsilon}.
$$

Let $\varepsilon_n$, $n\in \mathbb{N}$, be a sequence of positive numbers converging to $0$. For any $n\in\mathbb{N}$, let $F_n=F_{\varepsilon_n}$.
If
$\{u_n\}_{n=1}^{\infty}$ is a sequence contained in $\mathcal{A}_{\gamma}$ which
converges, as $n\to\infty$, to $u\in \mathcal{A}_{\gamma}$ in $L^1(B_R)$ and
satisfies $\lim_n F_n(u_n)=\lim_n\inf_{\mathcal{A}_{\gamma}} F_n$, then $u$ is a minimizer
for $F_0$ on $\mathcal{A}_{\gamma}$, that is $u$
solves the minimization problem \eqref{min0}.
\end{teo}

We conclude with the following remark. With the notation of Theorem~\ref{finalteo}, if
$\{u_n\}_{n=1}^{\infty}$ is a sequence contained in $\mathcal{A}_{\gamma}$ which
satisfies $\lim_n F_n(u_n)=\lim_n\inf_{\mathcal{A}_{\gamma}} F_n$, then, by 
Remark~\ref{compactnessoss}, we have that, up to passing to a subsequence,
$\{u_n\}_{n=1}^{\infty}$ actually
converges, as $n\to\infty$, to some function $u\in \mathcal{A}_{\gamma}$ in $L^1(B_R)$.




\section{Discussion}
We have provided a mathematical study of the inverse problem of photolithography.  The approach we
propose is to seek an approximate solution by formulating the geometrical problem using a
phase-field method.  We further relax the hard threshold involved in image exposure with an approximate
Heaviside function.  We show that the variational problem for the approximate solution is well-posed.
This opens a way into designing mathematically rigorous numerical methods.  We further show that
as the approximation parameter goes to zero, a theoretical limit, the original optimization problem
involving geometry is recovered.

\subsubsection*{Acknowledgements}
The authors learned about the inverse problem of photolithography from Apo Sezginer who gave
a seminar on this topic at the Institute for Mathematics and its Applications in 2004.  We
thank Dr.~Sezginer for helpful discussions.
Luca Rondi is partially supported by GNAMPA under 2008 and 2009 projects. Part of this work was done while Luca Rondi was visiting the School of Mathematics at the University of Minnesota, Minneapolis, USA, whose support and hospitality is gratefully
acknowledged.  Fadil Santosa's research is supported in part by NSF award DMS0807856.

\end{document}